\documentclass[11pt]{article}
\usepackage[english]{babel}
\usepackage{times}
\usepackage{paralist}
\usepackage{amsmath,amstext,amssymb,amsthm, amsfonts}
\usepackage{graphicx}
\usepackage{pstricks}

\newtheorem{lemma}{Lemma}[section]
\newtheorem{theorem}[lemma]{Theorem}
\newtheorem{corollary}[lemma]{Corollary}
\newtheorem{definition}[lemma]{Definition}

\newtheorem{remark}[lemma]{Remark}
\newtheorem{Assumption}[lemma]{Assumption}

\newcommand{\slim} {\mathop{\rm lim\,sup\,}}
\newcommand{\ilim} {\mathop{\rm lim\,inf\,}}

\def\U{\mathbb{U}}

\def\S{\mathbb{S}}
\def\K{\mathbb{K}}
\def\M{\mathbb{M}}
\def\Y{\mathbb{Y}}
\def\h{\mathbf{I}}
\def\X{\mathbb{X}}
\def\E{\mathbb{E}}
\def\A{\mathbb{A}}
\def\H{\mathbb{H}}
\def\R{\mathbb{R}}
\def\P{\mathbb{P}}
\def\F{\mathbb{F}}

\def\C{\mathbb{C}}
\def\W{\mathbb{W}}
\def\lll{\mathbb{L}}
\def\c{\bar{c}}
\def\B{\mathcal{B}}
\def\oo{\mathcal{O}}

\def\fff{\mathtt{F}}

\def\Psii{{\Xi_{\int}}}

\numberwithin{equation}{section}

\title{Markov Decision Processes with Incomplete Information and Semi-Uniform
Feller Transition Probabilities}
\begin{document}

\maketitle

\begin{center}
Eugene~A.~Feinberg \footnote{Department of Applied Mathematics and
Statistics,
 Stony Brook University,
Stony Brook, NY 11794-3600, USA, eugene.feinberg@sunysb.edu},\
Pavlo~O.~Kasyanov\footnote{Institute for Applied System Analysis,
National Technical University of Ukraine ``Igor Sikorsky Kyiv Polytechnic
Institute'', Peremogy ave., 37, build, 35, 03056, Kyiv, Ukraine,\
kasyanov@i.ua.},\ and Michael~Z.~Zgurovsky\footnote{National Technical University of Ukraine
``Igor Sikorsky Kyiv Polytechnic Institute'', Peremogy ave., 37, build, 1, 03056,
Kyiv, Ukraine,\
mzz@kpi.ua
}\\

\bigskip
\end{center}

\begin{abstract}
This paper deals with control of partially observable discrete-time stochastic systems.  It introduces and studies  Markov Decision Processes with Incomplete Information  and with semi-uniform Feller transition probabilities. The important feature  of these models is that their classic reduction to Completely Observable Markov Decision Processes with belief states preserves semi-uniform Feller continuity of transition probabilities. Under mild assumptions on cost functions,  optimal policies exist, optimality equations hold, and value iterations converge to optimal values for these  models. In particular, for Partially Observable Markov Decision Processes  the results of this paper imply new and generalize several known sufficient conditions on transition and observation probabilities for weak continuity of transition probabilities for Markov Decision Processes with belief states,  the existence of optimal policies, validity of optimality equations defining optimal policies, and convergence of value iterations to optimal values.
\end{abstract}

\section{Introduction}\label{S1} In many control problems the state of a controlled system is not known, and decision makers know only some information about the state. This takes place in many applications including signal processing, robotics, artificial intelligence, and medicine.  Except lucky exceptions, and Kalman's filtering is among them, problems with incomplete information are known to be difficult \cite{PT}.   The general approach to solving such problems was identified long ago in \cite{Ao, As, Dy, ShP}, and it is based on constructing a controlled system whose states are posterior state distributions for the original system.  These posterior distributions are often called belief probabilities or belief states.  Finding an optimal policy for a problem with incomplete state observation consists of two steps: (i) finding an optimal policy for the problem with belief states, and (ii) deriving from this policy an optimal policy for the original problem. This approach was introduced in \cite{Ao, As, Dy, ShP} for problems with finite state, observation, and action sets, and it holds for problems with Borel state, observation, and action sets \cite{Rh, Yu}.  If there is no optimal policy for the problem with belief states, then there is no optimal policy for the original problem.

  This paper deals with optimization of  expected total discounted costs for  discrete-time models.  We describe a large class of problems, for which optimal policies exist, satisfy optimality equations, which define optimal policies, and can be found by value iterations. In particular, this paper provides sufficient conditions for weak continuity of transition probabilities for models with belief states. For a particular model of Partially Observable Markov Decision  Process (POMDP), called ${\rm POMDP}_2$ in this paper, the related studies are \cite{FKZ,HL,Saldi,RS}. As known for long time, weak continuity of transition and observation probabilities for problems with incomplete information does not imply weak continuity of transition probabilities after the reduction to belief states. Examples are provided in \cite{FKZ}.

    Weak continuity of transition probabilities for models with belief states is an important property because these models are Markov Decision Processes (MDPs) with infinite state spaces.  Optimal policies  minimizing expected total discounted and undiscounted costs may not exist for such MDPs.  According to
  \cite[Theorem 2]{FKNMOR}, for MDPs with nonnegative costs and, if the discount factor is less than 1, with bounded below costs, weak continuity of transition probabilities and $\K$-inf-compactness of cost functions imply the existence of Markov optimal policies for finite-horizon problems and the existence of stationary optimal policies for infinite-horizon problems.  Under the mentioned two conditions, optimal policies satisfy optimality equations, and they can be found by value iteration starting from a zero value. For MDPs with belief states, $\K$-inf-compactness of cost functions follows from $\K$-inf-compactness of original cost functions \cite[Theorem 3.3]{FKZ}, and verifying  weak continuity of transition probabilities is a nontrivial matter.

\begin{figure}[h!]
  \center
  \includegraphics[width=110mm]{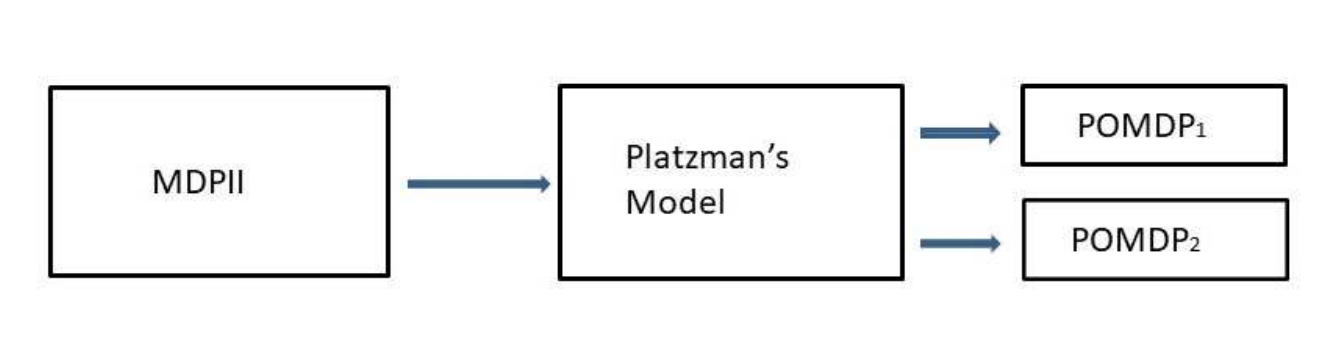}
  \caption{Relations between models of partially observable controlled Markov processes.  Platzman's model is defined as a particular case of an MDPII.  ${\rm POMDP}_1$ and ${\rm POMDP}_2$  are defined as particular cases of Platzman's model.}\label{Fig:Diagram}
\end{figure}

There are several models of controlled systems with incomplete state observations in the literature. Here we mostly consider a contemporary version of the original model introduced in  \cite{Ao, As, Dy, ShP} and called a Markov Decision Process with Incomplete Information (MDPII).  In this model the transitions are defined by transition probabilities $P(dw_{t+1},dy_{t+1}|w_t,y_t,a_t),$  where vectors $(w_t,y_t)$ represent states of the system at times $t=0,1,\ldots,$  $w_t$ and $y_t$ are unobservable and observable components of the state $(w_t,y_t)$, and $a_t$  are actions. In more contemporary studies the research focus switched to POMDPs. As was observed in \cite{Plat}, there are two different POMDP models in the literature, which we call ${\rm POMDP}_1$ and ${\rm POMDP}_2.$ For problems with finite state, observation, and control states, Platzman~\cite{Plat} introduced a ``plant'' model, which we adapt to problems with general state, observation, and control spaces and call Platzman's model. This model is more general than ${\rm POMDP}_1$ and ${\rm POMDP}_2;$ see Figure~\ref{Fig:Diagram}.

Platzman's model is a particular case of an MDPII when the transition probability does not depend on observations.
In other words, the transition probability in Platzman's model is $P(dw_{t+1},dy_{t+1}|w_t,a_t).$   ${\rm POMDP}_i,$  $i=1,2,$ are Platzman's models whose transition probabilities have special structural properties. These properties are $P(dw_{t+1},dy_{t+1}|w_t,a_t)=Q_1(dy_{t+1}|w_t,a_t)$ $P_1(dw_{t+1}| w_t,a_t)$ for ${\rm POMDP}_1$ and $P(dw_{t+1},dy_{t+1}|w_t,a_t)= Q_2(dy_{t+1}|a_t,w_{t+1}) P_2(dw_{t+1}| w_t,a_t)$ for ${\rm POMDP}_2,$ {where $P_i$ and $Q_i,$ $i=1,2,$ are transition and observation kernels respectively}. {Figure}~\ref{Fig:Diagram} illustrates the relations between  definitions of these four models based on the generality of the  transition probabilities $P(dw_{t+1},dy_{t+1}|w_t,y_t,a_t)$.  In particular, references \cite{Mo,SS,So} considered  ${\rm POMDP}_1,$ and references \cite{FKZ,HL,Saldi} considered ${\rm POMDP}_2.$

Belief-MDPs for MDPIIs are called Markov Decision Processes with Complete Information (MDPCIs) in this paper.  As mentioned above, the reduction of an MDPII  with Borel state, action, and observation sets to an MDPCI was introduced in \cite{Rh,Yu}. The reduction of a ${\rm POMDP}_2$ to a completely observable belief-MDP is described  in  \cite[Chapter 4]{HL}. The reduction of an MDPII to a ${\rm POMDP}_2$   described in \cite[Section 8.3] {FKZ} and the reduction of a ${\rm POMDP}_2$ to a completely observable belief-MDP described in  \cite[Chapter 4]{HL} also imply the reduction of an MDPII to an MDPCI.

This paper introduces the class of MDPIIs with semi-uniform Feller transition probabilities.  Theorem~\ref{th:mainMDPIINew} states that an MDPII has a transition probability from this class if and only if the transition probability of the corresponding MDPCI also belongs to this class. Theorem~\ref{th:mainMDPII} states similar results under more general conditions, which imply weaker continuity properties of value functions than the properties described in  Theorem~\ref{th:mainMDPIINew}. In view of Lemma~\ref{cor:Corollary 5.15.}, semi-uniform Feller transition probabilities are weakly continuous. In addition,  under mild conditions on cost functions described in Section~\ref{sec:MDPwithSemi-Feller}, there are optimal policies for MDPs with semi-uniform Feller transition probabilities.  This paper provides several sufficient conditions for the existence of optimal policies, validity of optimality equations, and convergence of value iterations.  In particular, the general theory implies the following sufficient conditions for weak continuity of transition probabilities for completely observable belief-MDPs corresponding to  POMDPs: (i) $P_i$ is weakly continuous and $Q_i,$  is continuous in total variation for  an $ {\rm POMDP}_i,$ $i=1,2$  (for $i=2$ this result was established in \cite{FKZ}); (ii) $P_2$ is continuous in total variation and $Q_2$ is continuous in total variation in the control parameter; sufficiency of continuity of $P_2$ in total variation was established in \cite{Saldi} for uncontrolled observation kernels, that is, $ Q_2(y_{t+1}|a_t,w_{t+1})= Q_2(y_{t+1}|w_{t+1})$.

Section~\ref{sec:MoDe} describes MDPIIs with  expected total costs, and Section~\ref{subsec:Reduction} describes their classic reduction to an MDPCI. Section~\ref{S4} introduces
semi-uniform Feller stochastic kernels {and it} provides the properties of semi-uniform Feller stochastic kernels.  In particular,
Lemma~\ref{cor:Corollary 5.15.} states  that semi-uniform Feller stochastic kernels are weakly continuous.  Semi-uniform Feller stochastic kernels were introduced  and studied
in \cite{FKZJThProb}, and some of the statements of Section~{\ref{S4}} are taken from there.
The basic known facts regarding the reduction of MDPIIs to MDPCIs are that this reduction preserves Borel measurability of transition probabilities \cite{Rh,Yu}, but it does
not preserve weak continuity of transition probabilities \cite[Examples 4.1 and 4.3]{FKZ}.  Section~\ref{sec:MDPwithSemi-Feller} describes the theory of MDPs
with the expected total costs and semi-uniform Feller transition probabilities.  { Theorem~\ref{prop:dcoe} establishes} the validity of optimality equations, convergence of value iterations {to optimal values},
existence of Markov optimal policies for finite horizon problems, and existence of stationary optimal  policies for infinite-horizon problems. {Related facts for MDPs with weakly and setwise continuous transition probabilities are \cite[Theorem 2]{FKNMOR} and \cite[Theorem 3.1]{FK} respectively.  MDPs with weakly and setwise continuous transition probabilities and with compact action sets were introduced and studied by Sch\"al \cite{sha, sha1, sha2}. Balder \cite{Ba} described a common approach to these models. MDPs with weakly and setwise continuous transition probabilities and possibly noncompact action sets were studied in \cite{FKNMOR} and \cite{FK, HL1} respectively. 
Weak continuity of transition probabilities is broadly used for problems with incomplete information, as described in this paper, and for inventory control~\cite{Finv}.} 
Section~\ref{sec:POMDPappl}
describes the results on the validity of optimality equations, convergence of value iterations {to optimal values}, and the existence of optimal policies for {belief-MDPs}  corresponding to MDPIIs,
Platzman's model, and POMDPs.  Proofs of several statements are presented in Appendix~\ref{app:B}.

Platzman's model in \cite{Plat}, references  \cite{FKZ, HL, SS, So} on POMDPs, and some papers on MDPIIs including \cite{Rh} considered one-step costs depending only on the unobservable states and actions.  References \cite{DY, FKZ,Yu} studied MDPIIs with one-step costs depending on unobservable states, observations, and actions.  In this paper we consider one-step costs   depending on unobservable states, observations, and actions.  Because of this, we consider in this paper more general POMDP models than are usually considered in the literature.  However, as shown in Section~\ref{sec:POMDPappl}, if one-step costs do not depend on observations, our results imply the known and new results for the classic Platzman's model~\cite{Plat} and POMDPs \cite{FKZ, HL, SS, So} with belief-MDPs having smaller state spaces $\P(\W)$ than state spaces $\P(\W{)}\times {\Y}$ for MDPCIs corresponding to Platzman's models, to POMDPs with one-step costs depending on observations, and to MDPIIs.  In general, costs may depend on observations in applications.  For example, for healthcare decisions during pandemics, costs depend not only on the health conditions of all the members {of} the population, which may be unknown, but also on the numbers of people with detected infections and {on} their conditions.

\section{Model Description}\label{sec:MoDe}

For a metric space $\S=(\S,\rho_\S),$ where $\rho_\S$ is a metric, let $\tau(\S)$ be the topology of $\S$ (the family of all open subsets of $\S$), and let ${\mathcal B}(\S)$ be its Borel
$\sigma$-field, that is, the $\sigma$-field generated by all open subsets of the
metric space $\S$. 
%
{For a subset $S$ of $\S$ let $\bar{S}$ denote the \textit{closure of} $S$ and $S^o$ the \textit{interior of} $S.$ Then $S^o\subset S\subset \bar{S},$ $S^o$ is open, and $\bar{S}$ is closed. Let $\partial S:=\bar{S}\setminus S^o$ denote the \textit{boundary of} $S.$} 
We denote by $\P(\S)$ the \textit{set of probability
measures} on $(\S,{\mathcal B}(\S)).$ A sequence of probability
measures $\{\mu^{\left(n\right)}\}_{n=1,2,\ldots}$ from $\P(\S)$
\textit{converges weakly} to $\mu\in\P(\S)$ if for every
bounded continuous function $f$ on $\S$
\[\int_\S f(s)\mu^{\left(n\right)}(ds)\to \int_\S f(s)\mu(ds) \qquad {\rm as \quad
}n\to\infty.
\]
A sequence of probability measures $\{\mu^{\left(n\right)}\}_{n=1,2,\ldots}$ from $\P(\S)$ \textit{converges in  total  variation} to $\mu\in\P(\S)$ if
\begin{equation}\label{eq:Kara1}
\begin{aligned}
\sup_{C\in \B(\S)}|\mu^{\left(n\right)}(C)-\mu(C)|\to 0\ {\rm as} \ n\to \infty;
\end{aligned}
\end{equation}
see~\cite{Steklov,UFL} for properties of these types of convergence of probability measures.
Note that $\P(\S)$ is a separable metric space with respect to the topology of weak convergence for probability measures, when $\S$ is
a separable metric space; \cite[Chapter~II]{Part}. Moreover,
according to Bogachev \cite[Theorem~8.3.2]{bogachev}, if the metric space $\S$ is separable, then the topology of weak convergence of probability measures on $(\S,\B(\S))$ coincides with the topology generated by the \textit{Kantorovich-Rubinshtein metric}
\begin{equation}\label{eq:KantorRubMetr}
\rho_{\P(\S)}(\mu,\nu):=\sup\left\{\int_\S f(s)\mu(ds)-\int_\S f(s)\nu(ds) \ \Big{|} \ f\in{\rm Lip}_1(\S),\ \sup_{s\in\S}|f(s)|\le 1 \right\},
\end{equation}
$\mu,\nu\in\P(\S),$ where
\[
{\rm Lip}_1(\S):=\{f:\S{\to}\mathbb{R}, \ |f(s_1)-f(s_2)|\le \rho_\S(s_1,s_2),\ \forall s_1,s_2\in\S\}.
\]

For a Borel subset $S$ of a metric space $(\S,\rho_\S),$ 
we always consider the
 metric space $(S,\rho_S),$ where $\rho_S:=\rho_\S\big|_{S\times S}.$  A subset $B$ of $S$ is called open
(closed) in $S$ if $B$ is open (closed respectively) in $(S,\rho_S)$. Of course, if $S=\S$, we omit ``in $\S$''. Observe that, in general, an open (closed) set
in $S$ may not be open (closed respectively). For $S\in\B(\S)$ we denote by $\B(S)$ the Borel $\sigma$-field on
$(S,\rho_S).$  Observe that $\B(S)=\{S\cap B:B\in\B(\S)\}.$

For metric spaces $\S_1$ and $\S_2$, a (Borel measurable) \textit{stochastic
kernel} $\Psi(ds_1|s_2)$ on $\S_1$ given $\S_2$ is a mapping $\Psi(\,\cdot\,|\,\cdot\,):\B(\S_1)\times \S_2{\to} [0,1]$, such that $\Psi(\,\cdot\,|s_2)$ is a
probability measure on $\S_1$ for any $s_2\in \S_2$, and $\Psi(B|\,\cdot\,)$ is a Borel measurable function on $\S_2$ for any Borel set $B\in\B(\S_1)$.  Another name for a stochastic kernel is a transition probability. A
stochastic kernel $\Psi(ds_1|s_2)$ on $\S_1$ given $\S_2$ defines a Borel measurable mapping $s_2\mapsto \Psi(\,\cdot\,|s_2)$ of $\S_2$ to the metric space
$\P(\S_1)$ endowed with the topology of weak convergence.
A stochastic kernel
$\Psi(ds_1|s_2)$ on $\S_1$ given $\S_2$ is called
\textit{weakly continuous (continuous in
total variation)}, if $\Psi(\,\cdot\,|s^{\left(n\right)})$ converges weakly (in
total  variation) to $\Psi(\,\cdot\,|s)$ whenever $s^{\left(n\right)}$ converges to $s$
in $\S_2$. For one-point sets $\{s_1\}\subset \S_1,$ we
sometimes write $\Psi(s_1|s_2)$ instead of $\Psi(\{s_1\}|s_2)$. Sometimes a weakly continuous stochastic kernel is called Feller, and a stochastic kernel continuous in total variation is called uniformly Feller \cite{Papa}.

Let $\S_1,\S_2,$ and $\S_3$ be Borel subsets of Polish spaces (a Polish space is a complete
separable metric space), and {let} $\Psi$ on $\S_1\times\S_2$ given $\S_3$ be a stochastic kernel. For each $A\in\B(\S_1),$
$B\in\B(\S_2),$ and $s_3\in\S_3,$ let
\begin{equation}\label{eq:marg_new}
\Psi(A,B|s_3):=\Psi(A\times B|s_3).
\end{equation}
In particular, we consider \textit{marginal} stochastic kernels
$\Psi(\S_1,\,\cdot\,|\,\cdot\,)$ on $\S_2$ given $\S_3$ and $\Psi(\,\cdot\,,\S_2|\,\cdot\,)$ on $\S_1$ given $\S_3.$

A \textit{Markov decision process with incomplete
information (MDPII)}  (Dynkin and  Yushkevich~\cite[Chapter~8]{DY},
Rhenius~\cite{Rh}, Yushkevich~\cite{Yu}; see also Rieder \cite{Ri}
and B\"auerle and Rieder~\cite{BR} for a version of this model
with transition probabilities having densities) is specified by a tuple
$(\W\times\Y,\A,P,c),$ where
\begin{itemize}
\item[(i)] $\W\times\Y$ is the \textit{state space}, where $\W$ and $\Y$ are Borel subsets of Polish spaces, and for $(w,y)\in\W\times\Y$ the unobservable component
    of the state $(w,y)$ is $w,$ and the observable component is $y;$
\item[(ii)] $\A$ is the \textit{action space}, which is assumed to be a Borel subset of a Polish space;
\item[(iii)] $P$ is a stochastic kernel on $\W\times\Y$ given
$\W\times\Y\times\A,$ which determines the distribution $P(\,\cdot\,|w,y,a)$ on $\W\times\Y$ of the new state, if $(w,y)\in\W\times\Y$ is the current state, and if $a\in A(y)$ is the current action, and it is assumed that the stochastic  kernel $P$ on {$\W\times\Y$} given $\W\times\Y\times\A$
is weakly continuous in $(w,y,a)\in \W\times\Y\times\A;$
\item[(iv)] $P_0(\,\cdot\,|w)$ is a stochastic kernel on $\Y$ given $\W,$ which determines the distribution of the observable part $y_0$ of the initial state, which may depend on the value of unobservable component $w_0=w$ of the initial state;
\item[(v)]
$c:\,\W\times\Y\times\A{\to}  \overline{\R}_+=[0,+\infty]$ is a Borel measurable \textit{one-step cost function}.
\end{itemize}

The Markov decision process with incomplete information \textit{evolves} as follows.
At time $t=0$, the unobservable component $w_0$ of the initial state has a given
prior distribution $p\in \P(\W).$ Let $y_0$ be the observable part
of the initial state. At
each time epoch $t=0,1,\ldots,$ if the state of the system is
$(w_t,y_t)\in\W\times\Y$ and the decision-maker chooses an action $a_t\in \A$,
then the cost $c(w_t,y_t,a_t)$ is incurred and
 the system moves to state $(w_{t+1},y_{t+1})$ according to
the transition law $P(\,\cdot\,|w_t,y_t,a_t).$

Define the \textit{observable histories}: $h_0:=y_0\in\H_0$
and
$h_t:=(y_0,a_0,y_1,a_1,\ldots,y_{t-1}, a_{t-1}, y_t)\in\H_t$ for all
$t=1,2,\dots,$
where $\H_0:=\Y$ and $\H_t:=\H_{t-1}\times \A\times \Y$
if $t=1,2,\dots$. Then a \textit{policy} for the MDPII is defined as
a sequence $\pi=\{\pi_t\}$ such that, for each $t=0,1,\dots,$
$\pi_t$ is a transition kernel on $\A$ given $\H_t$. Moreover, $\pi$
is called \textit{nonrandomized} if each probability measure
$\pi_t(\,\cdot\,|h_t)$ is concentrated at one point. The \textit{set of all policies} is denoted by $\Pi$.
The Ionescu Tulcea theorem (Bertsekas and Shreve \cite[pp.
140-141]{BS} or Hern\'andez-Lerma and Lasserre
\cite[p.178]{HLerma1}) implies that a policy $\pi\in \Pi,$
initial distribution $p\in \P(\W),$ initial state $y_0$
together with the transition
kernel $P$ determine a unique probability measure
$P_{p}^\pi$ on the set of all trajectories
$\mathbb{H}_{\infty}=(\W\times\Y\times \mathbb{A})^{\infty}$
endowed with the product $\sigma$-field defined by Borel
$\sigma$-fields of $\W$, $\Y$, and $\mathbb{A}$ respectively.
The expectation with respect to this probability measure is
denoted by $\E_{p}^\pi$.

Let us specify the performance criterion. For a finite horizon
$T=0,1,\ldots,$ and for a policy $\pi\in\Pi$, let  the
\textit{expected total discounted costs} be
\begin{equation}\label{eq1}
v_{T,\alpha}^{\pi}(p):=\mathbb{E}_{p}^{\pi}\sum\limits_{t=0}^{T-1}\alpha^tc(w_t,y_t,a_t),\quad
p\in \P(\W),
\end{equation}
where $\alpha\ge 0$ is the discount factor,
$v_{0,\alpha}^{\pi}(p)=0.$

When $T=\infty$, (\ref{eq1}) defines an
\textit{infinite horizon expected total discounted cost}, and we
denote it by $v_\alpha^\pi(p).$
For any function $g^{\pi}(p)$, including
$g^{\pi}(p)=v_{T,\alpha}^{\pi}(p)$ and
$g^{\pi}(p)=v_{\alpha}^{\pi}(p),$ define the \textit{optimal value}
$g(p):=\inf\limits_{\pi\in \Pi}g^{\pi}(p),$ $p\in\P(\W).$
For a given initial distribution $p\in\P(\W)$ of the initial unobservable component $w_0,$ a policy $\pi$ is called \textit{optimal} for the respective
criterion, if $g^{\pi}(p)=g(p)$ for all $p\in \P(\W).$ A policy is called
\emph{$T$-horizon discount-optimal} if $g^\pi=v_{T,\alpha}^\pi,$ and it
is called \emph{discount-optimal} if $g^\pi=v_{\alpha}^\pi.$

{We remark that the standard assumptions on the discount factor are either $\alpha\in [0,1)$ or $\alpha\in [0,1].$  However, since we assume that transition probabilities are weakly continuous and one-step costs are $\K$-inf-compact or satisfy a relaxed version of $\K$-inf-compactness stated in Definition~\ref{Ass:MeasKinf}, the same  monotonicity and continuity arguments apply to $\alpha>0;$ see the proof of Theorem 3 in \cite{FKNMOR}.  In addition, if    $\alpha\in [0,1),$ then it is possible to assume that $c$ is bounded from below rather than nonnegative.  This remark also applies for MDPs with setwise continuous transition probabilities $P$ and lower semi-continuous cost functions $c(x,a),$ which are inf-compact in variable $a;$ see \cite{FK}.  Of course, if $\alpha>1,$ then for many infinite-horizon problems the objective function is equal to $+\infty.$  The literature on MDPs with discount factors greater than 1 exists~\cite{HW}.  In particular, discount factors are relevant to opportunity costs and interest rates.  Discount factors greater than 1 are relevant to negative interest rates, which are offered by some banks at some countries.}

We recall that an MDP is defined by its state space, action space, transition probabilities, and one-step costs. An MDP is a particular case of an MDPII. Formally speaking, an MDP $(\X,\A,P,c)$ is an MDPII $(\W\times\Y,\A,P,c)$ with $\W$ being a singelton and $\Y=\X,$ where we follow the convention that $\W\times \X=\X$ in this case. In addition, for an MDP an initial state is observable. For an
MDP we consider an initial state $x$ instead of the initial { pair} $(P_0,p),$ where $p$ is the probability concentrated on a single point of which $\W$ consists. For an MDP, a nonrandomized policy is called \textit{Markov} if
all decisions depend only on the current state and time. A Markov
policy is called \textit{stationary} if all decisions depend only
on current states.

\section{Reduction of MDPIIs to MDPCIs}\label{subsec:Reduction}

In this section we formulate the well-known reduction of an MDPII $(\W\times\Y,\A,P,c)$ to a  belief-MDP (\cite{BS, DY, HLerma1, Rh, Yu}), which is called an MDPCI. For epoch $t=0,1,\ldots$ consider
the joint conditional probability $R(dw_{t+1}dy_{t+1}|z_t,y_t,a_t)$ on next state $(w_{t+1},y_{t+1})$ given the current state $(z_t,y_t)$ and the current control action $a_t$ defined by
\begin{equation}\label{3.3}
R(B\times C|z,y,a):=\int_{\W}P(B\times C|w,y,a)z(dw),
\end{equation}
$B\in \B(\W),$ $C\in\B(\Y),$ $(z,y,a)\in\P(\W)\times\Y\times\A.$ According to Bertsekas and Shreve \cite[Proposition~7.27]{BS}, there exists a stochastic kernel $H(z,y,a,y')[\,\cdot\,]=H(\,\cdot\,|z,y,a,y')$ on $\W$ given
$\P(\W)\times\Y\times \A\times\Y$ such that
\begin{equation}\label{3.4}
R(B\times C|z,y,a)=\int_CH(B|z,y,a,y')R(\W,dy'|z,y,a),
\end{equation}
$B\in \B(\W),$ $C\in\B(\Y),$ $(z,y,a)\in\P(\W)\times\Y\times\A.$ The stochastic kernel $H(\,\cdot\,|z,y,a,y')$ introduced in \eqref{3.4} defines a
measurable mapping $H:\,\P(\W)\times\Y\times \A\times \Y {\to}\P(\W).$ Moreover, the mapping $y'\mapsto H(z,y,a,y')$ is defined $R(\W,\,\cdot\,|z,y,a)$-a.s. uniquely for each triple $(z,y,a)\in \P(\W)\times\Y\times\A.$

Let $\h B$ denotes the \textit{indicator of an event} $B.$ The MDPCI is defined as an MDP with parameters
$(\P(\W)\times\Y,\A,q,\c),$ where
\begin{itemize}
\item[(i)] $\P(\W)\times\Y$ is the state space;
\item[(ii)] $\A$ is the
action set available at all states $(z,y)\in\P(\W)\times\Y;$
\item[(iii)] the
 one-step cost function $\c:\P(\W)\times\Y\times\A{\to}\overline{\R}$, defined as
\begin{equation}\label{eq:c}
\c(z,y,a):=\int_{\W}c(w,y,a)z(dw), \quad z\in\P(\W),\,y\in\Y,\, a\in\A;
\end{equation}
 \item[(iv)] $q$ on $\P(\W)\times\Y$
given $\P(\W)\times\Y\times \A$ is a stochastic kernel which determines the
distribution of the new
state as follows:
for $(z,y,a)\in \P(\W)\times\Y\times\A$
and for $D\in \B(\P(\W))$ and $C\in \B(\Y),$
\begin{equation}\label{3.7}
q(D\times C|z,y,a):=\int_{C}\h\{H(z,y,a,y')\in D\}
R(\W,dy'|z, y, a),
\end{equation}

\end{itemize}
see Yushkevich \cite{Yu}, Bertsekas and Shreve \cite[Corollary~7.27.1, p.~139]{BS}, or
Dynkin and Yushkevich \cite[p.~215]{DY} for details. Note that a particular  measurable
choice of a stochastic kernel $H$ from (\ref{3.4}) does not { effect}
the definition of $q$ in (\ref{3.7}).

There is a correspondence between the policies for an  MDPII $(\W\times\Y,\A,P,c)$  and for the corresponding MDPCI $(\P(\W)\times\Y,\A,q,\c)$ in the sense that for a policy in one of these models there exists a policy in another model with the same expected total costs; see \cite{Rh,Yu} or \cite[Section 4.3]{HL}. In Section~\ref{sec:POMDPappl} we provide sufficient conditions for the existence of an optimal policy in the MDPCI $(\P(\W)\times\Y,\A,q,\c)$ in terms of the assumptions on the initial MDPII $(\W\times\Y,\A,P,c)$ and apply the results to Platzman's model and POMDPs.  In particular, under natural conditions the existence of optimal policies and validity of optimality equations and value iterations for MDPCIs follow from Theorem~\ref{prop:dcoe}. {For problems with finite and infinite horizons,} if $\phi$ is a Markov optimal policy for the MDPCI, then an optimal  policy $\pi$ for the MDPII can be defined as
$a_t=\pi_t(h_t)=\phi_t(z_t,y_t),$ where $z_t$ is the posterior distribution of the unobservable component $w_t$ of the state $x_t$ given the observations $h_t=(y_0,a_0,\ldots,y_{t-1},a_{t-1},y_t),$ the initial distribution $p$ of $w_0,$ and  $t>0.$  As discussed in Section~\ref{sec:POMDPappl}, for Paltzman's models and, in particular, for POMDPs, the values of $\phi_t(z_t,y_t)$ can be selected independent of $y_t$ if one-step costs do not depend on observations. {For infinite-horizon MDPs usually there exist stationary optimal policies, and the described scheme applies to them since stationary policies are Markov.}


\section{Semi-Uniform
Feller Stochastic Kernels and their {Properties}}\label{S4}

In this section we formulate the semi-uniform Feller property for stochastic kernels {and} describe {its} basic properties. In particular, Theorem~\ref{th:equivWTV} provides {its} equivalent definitions. Theorem~\ref{th:concept} establishes a necessary and sufficient condition for a stochastic kernel to be {semi-uniform Feller}. This condition is Assumption~\ref{AssKern}, whose stronger version was introduced in \cite[Theorem~4.4]{Steklov}. Theorem~\ref{th:extra} describes the preservation of {semi-uniform Fellerness} under the integration operation.


Let $\S_1,$ $\S_2,$ and $\S_3$ be Borel subsets of Polish spaces, and $\Psi$ on $\S_1\times\S_2$ given $\S_3$ be a stochastic kernel.
\begin{definition}\label{defi:unifFP}{{\rm(Feinberg et al. \cite{FKZJThProb})}}
A stochastic kernel $\Psi$ on $\S_1\times\S_2$ given $\S_3$ is semi-uniform Feller if, for each sequence $\{s_3^{\left(n\right)}\}_{n=1,2,\ldots}\subset\S_3$ that converges to $s_3$ in $\S_3$ and for each bounded continuous function $f$ on $\S_1,$
\begin{equation}\label{eq:equivWTV3}
\lim_{n\to\infty} \sup_{B\in \B(\S_2)} \left| \int_{\S_1} f(s_1) \Psi(ds_1,B|s_3^{\left(n\right)})-\int_{\S_1} f(s_1) \Psi(ds_1,B|s_3)\right|= 0.
\end{equation}
\end{definition}

We recall that the marginal measure $\Psi(ds_1,B|s_3),$ $s_3\in\S_3,$ is defined in \eqref{eq:marg_new}.
The term ``semi-uniform'' is used in Definition~\ref{defi:unifFP} because the uniform property holds in \eqref{eq:equivWTV3} only with respect to the {second} coordinate. If the uniform property holds with respect to both  coordinates, then the stochastic kernel $\Psi$ on $\S_1\times\S_2$ given $\S_3$ is continuous in total variation, and it is sometimes called uniformly Feller \cite{Papa}. 

{\begin{lemma}\label{cor:Corollary 5.15.} A semi-uniform Feller stochastic kernel $\Psi$ on $\S_1\times \S_2$ given $\S_3$ is weakly continuous.
\end{lemma}
\begin{proof}
Definition~\ref{defi:unifFP} implies that for each sequence $\{s_3^{\left(n\right)}\}_{n=1,2,\ldots}\subset\S_3$ that converges to $s_3$ in $\S_3,$ for each bounded continuous function $f$ on $\S_1,$ and for each $B\in \B(\S_2)$
\[\lim_{n\to\infty}
 \int_{\S_1} f(s_1) \Psi(ds_1,B|s_3^{\left(n\right)})=\int_{\S_1} f(s_1) \Psi(ds_1,B|s_3),
\]
and, in view of Sch\"al \cite[Theorem~3.7(iii,viii)]{sha}, this property implies weak continuity of $\Psi$ on $\S_1\times \S_2$ given $\S_3.$
\end{proof}
}

Let us consider some basic definitions.
\begin{definition}\label{defi:semi}
Let $\S$ be a metric space.
 A function $f:\S{\to} \mathbb{R}$ is called
\begin{itemize}
\item[{\rm(i)}] \textit{lower semi-continuous} (l.s.c.) at a point $s\in\S$ if $\mathop{\ilim}\limits_{s'\to s}f(s')\ge f(s);$
\item[{\rm(ii)}] \textit{upper semi-continuous} at $s\in\S$ if $-f$ is lower semi-continuous at $s;$
\item[{\rm(iii)}]
\textit{continuous} at $s\in\S$ if $f$ is both lower and upper semi-continuous at $s;$
\item[{\rm(iv)}] \textit{lower / upper semi-continuous (continuous respectively) (on $\S$)} if $f$ is lower / upper semi-continuous (continuous respectively) at each $s\in\S.$
\end{itemize}
\end{definition}
For a metric space $\S,$ let  $\F(\S),$ $\lll(\S),$ and $\C(\S)$ be the spaces of all
real-valued functions, all real-valued lower semi-continuous functions, and all real-valued continuous functions respectively defined
on the metric space $\S.$ The following definitions are taken from \cite{FKL2}.

\begin{definition}\label{defi:equi}
A {family} $\mathtt{F}\subset \F(\S)$ of real-valued functions on a metric space $\S$ is called
\begin{itemize}
\item[{\rm(i)}]
\textit{lower semi-equicontinuous at a point} $s\in \S$ if $\ilim_{s'\to s}\inf_{f\in\mathtt{F}} (f(s')-f(s))\ge0;$
\item[{\rm(ii)}]
\textit{upper semi-equicontinuous at a point} $s\in \S$ if the {family} $\{-f\,:\,f\in\mathtt{F}\}$ is lower semi-equicontinuous at $s\in \S;$
\item[{\rm(iii)}] \textit{equicontinuous at a point $s\in\S$}, if $\mathtt{F}$ is both lower and upper semi-equicontinuous at $s\in\S,$ that is, $\mathop{\lim}\limits_{s'\to s} \mathop{\sup}\limits_{f\in\mathtt{F}} |f(s')-f(s)|=0;$
\item[{\rm(iv)}]
\textit{lower / upper semi-equicontinuous (equicontinuous respectively)} (\textit{on $\S$}) if it is lower / upper semi-equicontinuous (equicontinuous respectively) at all $s \in \S;$
\item[{\rm(v)}] \textit{uniformly bounded (on $\S$)}, if there exists a constant $M<+\infty $ such that $ |f(s)|\le M$ for all $s\in\S$ and  for
all $f\in \mathtt{F}.$
\end{itemize}
\end{definition}

Obviously, if a {family}
$\mathtt{F}\subset \F(\S)$ is lower semi-equicontinuous, then
$\mathtt{F}\subset \lll(\S).$ Moreover, if a {family} $\mathtt{F}\subset \F(\S)$ is equicontinuous, then $\mathtt{F}\subset \C(\S).$

\subsection{{Basic Properties of Semi-Uniform Feller} Stochastic Kernels}\label{subsec:KernelsInProduct}

Let $\S_1$, $\S_2,$ and $\S_3$ be Borel subsets of Polish spaces, and let $\Psi$ on $\S_1\times\S_2$ given $\S_3$ be a stochastic kernel.
For each set $A\in\B(\S_1)$ consider the {family} of functions
\begin{equation}\label{eq:familyoffunctions}
\fff^\Psi_A=\{  s_3\mapsto \Psi(A\times B |s_3):\, B\in \B(\S_2)\}
\end{equation}
mapping $\S_3$ into $[0,1].$ Consider the following type of continuity for stochastic kernels on $\S_1\times\S_2$ given $\S_3.$
\begin{definition}\label{defi:wtv}
A stochastic kernel $\Psi$ on $\S_1\times\S_2$ given $\S_3$ is called \textit{WTV-continuous}, if for each $\oo \in\tau (\S_1)$
the {family} of functions $\fff^\Psi_\oo$ is lower semi-equicontinuous on $\S_3.$
\end{definition}
Definition~\ref{defi:equi} directly implies that the stochastic kernel $\Psi$ on $\S_1\times\S_2$ given $\S_3$ is WTV-continuous if and only if for each $\oo \in\tau (\S_1)$
\begin{equation}\label{eq:equivWTV0}
\ilim_{n\to\infty} \inf_{B\in \B(\S_2)\setminus\{\emptyset\}} \left( \Psi(\oo \times B|s_3^{\left(n\right)})-\Psi(\oo \times B|s_3)\right)\ge 0,
\end{equation}
whenever $s_3^{\left(n\right)}$ converges to $s_3$ in $\S_3.$

{Since} $\emptyset\in\B(\S_2),$ \eqref{eq:equivWTV0} holds if and only if
\begin{equation}\label{eq:wtvsc}
\lim_{n\to\infty} \inf_{B\in \B(\S_2)} \left( \Psi(\oo \times B|s_3^{\left(n\right)})-\Psi(\oo \times B|s_3)\right)= 0.
\end{equation}
{WTV-continuity of the stochastic kernel
$\Psi$ on $\S_1\times\S_2$ given $\S_3$ implies  continuity in total variation of its marginal kernel $\Psi(\S_1,\,\cdot\,|\,\cdot\,)$ on $\S_2$ given $\S_3$
because
\begin{equation*}
\lim_{n\to\infty} \sup_{B\in \B(\S_2)} \left| \Psi(\S_1 \times B|s_3^{\left(n\right)})-\Psi(\S_1 \times B|s_3)\right|=\lim_{n\to\infty} \sup_{B\in \B(\S_2)} \left( \Psi(\S_1 \times B|s_3^{\left(n\right)})-\Psi(\S_1 \times B|s_3)\right)= 0,
\end{equation*}
where the second equality
follows from equality \eqref{eq:wtvsc} with $\oo:=\S_1$ and from $\Psi(\S_1\times\S_2|\,\cdot\,)=1.$
}


Similarly to Parthasarathy \cite[Theorem~II.6.1]{Part}, where the necessary and sufficient conditions for weakly convergent probability
measures were considered, the following theorem
provides several useful equivalent definitions of the {semi-uniform Feller} stochastic kernels.

\begin{theorem}{\rm(Feinberg et al \cite[Theorem~3]{FKZJThProb})}\label{th:equivWTV}
For a stochastic kernel $\Psi$ on $\S_1\times\S_2$ given $\S_3$ the following conditions are equivalent:
\begin{itemize}
\item[{\rm(a)}] the stochastic kernel $\Psi$ on $\S_1\times\S_2$ given $\S_3$ is semi-uniform Feller;
\item[{\rm(b)}] the stochastic kernel $\Psi$ on $\S_1\times\S_2$ given $\S_3$ is WTV-continuous;
\item[{\rm(c)}] if $s_3^{\left(n\right)}$ converges to $s_3$ in $\S_3,$ then for each closed set $C$ in $\S_1$
\begin{equation}\label{eq:wtvscConv}
\lim_{n\to\infty} \sup_{B\in \B(\S_2)} \left( \Psi(C \times B|s_3^{\left(n\right)})-\Psi(C \times B|s_3)\right)= 0;
\end{equation}
\item[{\rm(d)}] if $s_3^{\left(n\right)}$ converges to $s_3$ in $\S_3,$ then, for each $A\in\B(\S_1)$ such that $\Psi(\partial A,\S_2|s_3)=0,$
\begin{equation}\label{eq:equivWTV2}
\lim_{n\to\infty} \sup_{B\in \B(\S_2)} | \Psi(A \times B|s_3^{\left(n\right)})-\Psi(A \times B|s_3)|= 0;
\end{equation}
\item[{\rm(e)}] if $s_3^{\left(n\right)}$ converges to $s_3$ in $\S_3,$ then, for each nonnegative bounded lower semi-continuous function $f$ on $\S_1,$
\begin{equation}\label{eq:equivWTV4}
\ilim_{n\to\infty} \inf_{B\in \B(\S_2)} \left( \int_{\S_1} f(s_1) \Psi(ds_1,B|s_3^{\left(n\right)})-\int_{\S_1} f(s_1) \Psi(ds_1,B|s_3)\right)= 0;
\end{equation}
\end{itemize}
{and each of these conditions implies continuity in total variation of the marginal kernel $\Psi(\S_1,\,\cdot\,|\,\cdot\,)$ on $\S_2$ given $\S_3.$}
\end{theorem}

Note that, since $\emptyset\in\B(\S_2),$ \eqref{eq:wtvscConv} holds if and only if
\begin{equation}\label{eq:equivWTV1}
\slim_{n\to\infty} \sup_{B\in \B(\S_2)\setminus\{\emptyset\}} \left( \Psi(C \times B|s_3^{\left(n\right)})-\Psi(C \times B|s_3)\right)\le 0,
\end{equation}
and similar remarks are applicable to \eqref{eq:equivWTV2} and \eqref{eq:equivWTV4} with the inequality ``$\ge$'' taking place in \eqref{eq:equivWTV4}.

%

Let us consider the following assumption. According to Feinberg et al \cite[Example~1]{FKZJThProb}, Assumption~\ref{AssKern} is weaker than combined assumptions (i) and (ii) in \cite[Theorem~4.4]{Steklov}, where the base $\tau_b^{s_3}(\S_1)$ is the same for all $s_3\in\S_3.$

\begin{Assumption}\label{AssKern} Let $\Psi$ be a stochastic kernel on $\S_1\times\S_2$ given $\S_3,$ and
let for each $s_3\in\S_3$ the topology on $\S_1$ have a countable base
$\tau_b^{s_3}(\S_1)$ such that:
\begin{itemize}
 \item[{\rm(i)}] $\S_1\in\tau_b^{s_3}(\S_1);$
 \item[{\rm(ii)}]  for
each finite intersection $\oo=\cap_{i=1}^ k {\oo}_{i},$ $k=1,2,\ldots,$ of sets
$\oo_{i}\in\tau_b^{s_3}(\S_1),$ $i=1,2,\ldots,k,$
the family of functions $\fff^\Psi_\oo,$ defined in \eqref{eq:familyoffunctions}, is equicontinuous at $s_3.$
\end{itemize}
\end{Assumption}

Note that Assumption~\ref{AssKern}(ii) holds if and only if for
each finite intersection $\oo=\cap_{i=1}^ k {\oo}_{i}$ of sets
$\oo_{i}\in\tau_b^{s_3}(\S_1),$ $i=1,2,\ldots,k,$
\begin{equation}\label{eq:equivWTV0new1}
\lim_{n\to\infty} \sup_{B\in \B(\S_2)} \left| \Psi(\oo \times B|s_3^{\left(n\right)})-\Psi(\oo \times B|s_3)\right|= 0
\end{equation}
if $s_3^{\left(n\right)}$ converges to $s_3$ in $\S_3.$

Theorem~\ref{th:concept} shows that Assumptions~\ref{AssKern} is a necessary and sufficient condition for semi-uniform Feller continuity.

\begin{theorem}{\rm(Feinberg et al \cite[Theorem~4]{FKZJThProb})}\label{th:concept}
The stochastic kernel $\Psi$ on $\S_1\times\S_2$ given $\S_3$ is semi-uniform Feller if and only if it satisfies Assumption~\ref{AssKern}.
\end{theorem}

Now let $\S_4$ be a Borel subset of a Polish space, and let $\Xi$ be a stochastic kernel on $\S_1\times\S_2$ given $\S_3\times\S_4.$ Consider the stochastic kernel $\Psii$ on $\S_1\times\S_2$ given $\P(\S_3)\times\S_4$ defined by
\begin{equation}\label{eq:extra1}
\Psii(A\times B|\mu,s_4):=\int_{\S_3}\Xi(A\times B |s_3,s_4)\mu(ds_3),\quad A\in\B(\S_1),\,B\in\B(\S_2),\,\mu\in\P(\S_3),\,s_4\in\S_4.
\end{equation}
{We observe that \eqref{eq:extra1} becomes \eqref{3.3} with $\Psii:= R,$ $\Xi:=P,$ $\S_1:=\W,$ $\S_2:=\Y,$ $\S_3:=\W,$ and $\S_4:=\Y\times\A.$ This is our
main motivation for writing \eqref{eq:extra1}.}

The following theorem establishes the preservation of {semi-uniform Fellerness} of the integration operation in \eqref{eq:extra1}.

\begin{theorem}{\rm(Feinberg et al \cite[Theorem~5]{FKZJThProb})}\label{th:extra}
The stochastic kernel $\Psii$ on
$\S_1\times\S_2$ given $\P(\S_3)\times\S_4$ is {semi-uniform Feller}
if and only if \, $\Xi$ on $\S_1\times\S_2$ given $\S_3\times\S_4$ is {semi-uniform Feller}.
\end{theorem}
\subsection{Continuity Properties of Posterior Distributions}\label{subsec:ContBR}

In this subsection we describe sufficient conditions for semi-uniform Feller continuity of posterior distributions. The main result of this section is Theorem~\ref{th:CBRmain}.

Let $\S_1,$ $\S_2,$ and $\S_3$ be Borel subsets of Polish spaces, and $\Psi$ on $\S_1\times\S_2$ given $\S_3$ be a stochastic kernel. By Bertsekas and Shreve \cite[Proposition~7.27]{BS}, there exists a stochastic kernel $\Phi$ on $\S_1$ given
$\S_2\times\S_3$ such that
\begin{equation}\label{eq:CBR1}
\Psi(A\times B|s_3)=\int_B\Phi(A|s_2,s_3)\Psi(\S_1,ds_2|s_3),\quad A\in \mathcal{B}(\S_1),\  B\in \mathcal{B}(\S_2),\ s_3\in\S_3.
\end{equation}

The stochastic kernel $\Phi(\,\cdot\,|s_2,s_3)$ on $\S_1$ given
$\S_2\times\S_3$ defines a measurable
mapping $\Phi:\,\S_2\times\S_3 \to\P(\S_1),$ where
$\Phi(s_2,s_3)(\,\cdot\,)=\Phi(\,\cdot\,|s_2,s_3).$ According to Bertsekas and Shreve \cite[Corollary~7.27.1]{BS}, for each $s_3\in
\S_3$ the mapping $\Phi(\,\cdot\,,s_3):\S_2\to\P(\S_1)$ is defined
$\Psi(\S_1,\,\cdot\,|s_3)$-almost surely uniquely in $s_2\in\S_2.$
Let us consider the stochastic kernel  $\phi$ defined by
\begin{equation}\label{eq:CBR2}
\phi(D\times B|s_3):=\int_{B}\h\{\Phi(s_2,s_3)\in D\}\Psi(\S_1,ds_2|s_3),\quad D\in \mathcal{B}(\P(\S_1)),\ B\in\B(\S_2),\ s_3\in\S_3,
\end{equation}
 where a particular choice of
a stochastic kernel $\Phi$ satisfying (\ref{eq:CBR1}) does not effect the
definition of $\phi$ in (\ref{eq:CBR2}).

In models for decision making with incomplete information, $\phi$ is the transition probability between belief states, which are posterior distributions of states; \eqref{3.7}. Continuity properties of $\phi$ play the fundamental role in the studies of models with incomplete information. Theorem~\ref{th:CBRmain} characterizes such properties, and this is the reason for the title of this section. Let us consider the following assumption.
\begin{Assumption}\label{Ass:H}  {For a stochastic kernel $\Psi$ on $\S_1\times\S_2$ given $\S_3,$ there} exists a stochastic
kernel $\Phi$ on $\S_1$ given $\S_2\times\S_3$ satisfying
(\ref{eq:CBR1}) such that, if a sequence
$\{s_3^{\left(n\right)}\}_{n=1,2,\ldots}\subset\S_3$ converges to
$s_3\in\S_3$ as $n\to\infty,$ then there exists a
subsequence $\{s_3^{\left(n_k\right)}\}_{k=1,2,\ldots}\subset
\{s_3^{\left(n\right)}\}_{n=1,2,\ldots}$ and a measurable subset $B$ of
$\,\S_2$ such that
\begin{equation}\label{eq:CBR3}
\Psi(\S_1\times B|s_3)=1\quad\mbox{and}\quad\Phi(s_2,s_3^{\left(n_k\right)})\mbox{ converges weakly to }\Phi(s_2,s_3)\quad\mbox{for all }s_2\in B.
\end{equation}
In other words, the convergence in \eqref{eq:CBR3} holds $\Psi(\S_1,\,\cdot\,|s_3)$-almost
surely.
\end{Assumption}

{According to Theorem~9.2.1 from \cite{Dud} {stating} the relation between convergence in probability and almost sure convergence, Assumption~\ref{Ass:H} holds if and only if the following statement holds: if a sequence
$\{s_3^{\left(n\right)}\}_{n=1,2,\ldots}\subset\S_3$ converges to
$s_3\in\S_3$ as $n\to\infty,$ then
\begin{equation}\label{eq:CBR4}
\rho_{\P(\S_1)}(\Phi(s_2,s_3^{\left(n\right)}),\Phi(s_2,s_3))\to0\mbox{ in probability }\Psi(\S_1,ds_2|s_3),
\end{equation}
where $\rho_{\P(\S_1)}$ is an arbitrary metric that induces the topology of weak convergence of probability measures on $\S_1,$ and, in particular, $\rho_{\P(\S_1)}$ can be the Kantorovich-Rubinshtein metric defined in \eqref{eq:KantorRubMetr}.}

The following theorem, which is the main result of this section, provides necessary and sufficient conditions for {semi-uniform Fellerness} of a stochastic kernel $\phi$ in terms of the properties of a given stochastic kernel $\Psi.$ This theorem and the results of Subsection~\ref{subsec:KernelsInProduct} provide the necessary and sufficient conditions for the semi-uniform Feller property of the MDPCIs in terms of the conditions on the transition kernel in the initial model for decision making with incomplete information.

\begin{theorem}\label{th:CBRmain}
For a stochastic kernel $\Psi$ on $\S_1\times\S_2$ given $\S_3$ the following conditions are equivalent:
\begin{itemize}
\item[{\rm(a)}] the stochastic kernel $\Psi$ on $\S_1\times\S_2$ given $\S_3$ is semi-uniform Feller;
\item[{\rm(b)}] {the marginal kernel $\Psi(\S_1,\,\cdot\,|\,\cdot\,)$ on $\S_2$ given $\S_3$ is continuous in total variation and} Assumption~\ref{Ass:H} holds;
\item[{\rm({c})}] the stochastic kernel $\phi$ on $\P(\S_1)\times \S_2$ given $\S_3$ is semi-uniform {Feller.} 
\end{itemize}
\end{theorem}
\begin{proof}
See Appendix \ref{app:B}.
\end{proof}

\section{Markov Decision Processes with Semi-Uniform Feller Kernels
}\label{sec:MDPwithSemi-Feller}

Let $\X_W$ and $\X_Y$ be Borel subsets of Polish spaces. In this section we consider the special class of MDPs with semi-uniform Feller transition kernels, when {the state space is} $\X:=\X_W\times\X_Y.$ These results are important for MDPIIs with semi-uniform Feller transition kernels from Section~\ref{sec:POMDPappl}, where $\X_W:=\P(\W)$ and $\X_Y=\Y.$

For an $\overline{\mathbb{R}}$-valued
function $f,$ defined on a nonempty subset $U$ of a metric
space $\mathbb{U},$ consider the level sets
\begin{equation}\label{def-D}
\mathcal{D}_f(\lambda;U)=\{y\in U \, : \,  f(y)\le
\lambda\},\qquad \lambda\in\R.
\end{equation}  We recall that a
function $f$ is \textit{inf-compact on $U$} if all the
level sets $\mathcal{D}_f(\lambda;U)$
are compact.

For a metric space $\U$, we denote by $\K(\U)$  the family of all
nonempty compact subsets of $\U.$
\begin{definition}{\rm(Feinberg et al. \cite[Definition~1.1]{FKN})}\label{K-inf-comp}
A function $u:\S_1\times\S_2{\to} \overline{\mathbb{R}}$ is called $\K$-inf-compact if this function is inf-compact on $K \times\S_2$ for each $K\in \K(\S_1).$
\end{definition}

The fundamental importance of $\K$-inf-compactness is that Berge's theorem stating lower semicontinuity of the value function holds for possibly noncompact action sets; Feinberg et al \cite[Theorem~1.2]{FKN}. In particular, {this fact} allows us to
consider the MDPII $(\W\times\Y,\A,P,c)$ with a possibly noncompact action space $\A$ and unbounded one-step cost $c$ and examine convergence of value  iterations for this model in Theorem~\ref{th:mainMDPII}, for Platzman's model in Corollaries~\ref{cor:Plpl}, \ref{cor:POMDPmaincwa}, and for POMDPs in Corollaries~\ref{cor:POMDPmain1}, \ref{cor:POMDPmain}.

\begin{definition}\label{Ass:MeasKinf}
A Borel measurable function $u:\S_1\times\S_2\times\S_3{\to} \overline{\mathbb{R}}$ is called {measurable $\K$-inf-compact on ${(}\S_1\times\S_{3},\S_{2}{)}$ or $\M\K(\S_1\times\S_3,\S_2)$-inf-compact } if for each $s_2\in \S_2$ the function $(s_1,s_3)\mapsto u(s_1,s_2{,}s_3)$ is $\K$-inf-compact on $\S_1\times\S_3.$
\end{definition}


Consider a discrete-time MDP $(\X,\A,q,c)$ with a state space $\X=\X_W\times \X_Y,$ an
action space $\mathbb{A},$ one-step costs $c,$ and
transition probabilities $q.$ Assume that $\X_W,\X_Y,$ and $\mathbb{A}$ are
\textit{Borel subsets} of Polish spaces. Let $LW(\mathbb{X})$ be the class of all nonnegative  Borel measurable
functions $\varphi:\mathbb{X}\to\overline{\mathbb{R}}$
such that $w\mapsto \varphi(w,y)$ is lower semi-continuous on $\X_W$ for each $y\in\X_Y.$ For any $\alpha\ge 0$ and $u\in LW(\X),$ we consider
\begin{equation}\label{e:defeta}
\eta_u^\alpha(x,a)=c(x,a)+\alpha\int_\X  u(\tilde{x})q(d\tilde{x}| x,a),
\quad(x,a)\in \X\times\A.
\end{equation}

The following theorem is the main result of this
section. It states the validity of optimality equations, convergence of value iterations, and existence of optimal policies for MDPs with semi-uniform Feller transition
probabilities and  {$\M\K(\W\times\A,\Y)$}-inf-compact one-step cost functions, when the goal is to minimize expected total costs.  For MDPs with weakly continuous
transition probabilities  the similar result is \cite[Theorem 2]{FKNMOR}, and for MDPs with setwise continuous transition probabilities the similar result is  \cite[Theorem 3.1]{FK}.
 {Theorem~\ref{prop:dcoe} does not follow from these two results.  In particular, the cost function is lower semi-continuous in \cite[Theorem 2]{FKNMOR}.
 The corresponding assumption for Theorem~~\ref{prop:dcoe} would be lower semi-continuity of the cost function $c,$  but the function $c(w,y,a)$ may not be lower semi-continuous in $y.$
 \cite[Theorem 3.1]{FK} assumes setwise continuity of the transition probability $q$ in the control parameter, which may not hold in this paper.}
Theorem~\ref{prop:dcoe} is applied in Theorem~\ref{th:mainMDPII} to MDPCIs $(\P(\W)\times\Y,\A,q,\c)$.
\begin{theorem}\label{prop:dcoe}{\rm(Expected Total Discounted  Costs)} Let us consider an MDP $(\X,\A,q,c)$ with
$\X=\X_W\times \X_Y,$ for each $y\in\X_Y$ the stochastic kernel $q(\,\cdot\,|\,\cdot\,,y,\,\cdot\,)$ on $\X$ given $\X_W\times\A$ being semi-uniform Feller, and the nonnegative function $c:\X\times \A{\to} \overline{\mathbb{R}}$ being
 {$\M\K(\X_W\times\A,\X_Y)$}-inf-compact.
Then
\begin{itemize}
\item[{\rm(i)}] the functions $v_{t,\alpha},$ $t=0,1,\ldots,$ and $v_\alpha$
{belong} to $LW(\mathbb{X}),$ and $v_{t,\alpha}(x)\uparrow
v_\alpha (x)$ as $t \to +\infty$ for all $x\in \X;$
\item[{\rm(ii)}] $v_{t+1,\alpha}(x)=\min\limits_{a\in \A}\eta_{v_{t,\alpha}}^\alpha(x,a),$
$x\in \X,$ $t=0,1,...,$ where $v_{0,\alpha}(x)=0$ for all $x\in \X,$ and the nonempty sets $A_{t,\alpha}(x):=\{a\in \A\,:\,v_{t+1,\alpha}(x)=\eta_{v_{t,\alpha}}^\alpha(x,a) \},$ $x\in \X,$ $t=0,1,\ldots,$ satisfy the following properties: (a) the graph
${\rm Gr}(A_{t,\alpha})=\{(x,a)\,:\, x\in\X,\, a\in A_{t,\alpha}(x)\},$
$t=0,1,\ldots,$ is a Borel subset of $\X\times \mathbb{A},$ and (b)
if $v_{t+1,\alpha}(x)=+\infty,$ then $A_{t,\alpha}(x)=\A$ and, if
$v_{t+1,\alpha}(x)<+\infty,$ then $A_{t,\alpha}(x)$ is compact;
\item[{\rm(iii)}]  for any $T=1,2,\ldots,$ there exists a Markov optimal
$T$-horizon policy $(\phi_0,\ldots,\phi_{T-1}),$ and, if for an
$T$-horizon Markov policy $(\phi_0,\ldots,\phi_{T-1})$ the
inclusions $\phi_{T-1-t}(x)\in A_{t,\alpha}(x),$ $x\in\X,$
$t=0,\ldots,T-1,$ hold, then this policy is $T$-horizon optimal;
\item[{\rm(iv)}] $v_{\alpha}(x)=\min\limits_{a\in \A}\eta_{v_{\alpha}}^\alpha(x,a),$ $x\in \X,$
and the nonempty sets $A_{\alpha}(x):=\{a\in
\A\,:\,v_{\alpha}(x)=\eta_{v_{\alpha}}^\alpha(x,a) \},$ $x\in \X,$
satisfy the following properties: (a) the graph ${\rm
Gr}(A_{\alpha})=\{(x,a)\,:\, x\in\X,\, a\in A_\alpha(x)\}$  is a Borel
subset of $\X\times \mathbb{A},$ and (b) if $v_{\alpha}(x)=+\infty,$
then $A_{\alpha}(x)=\A$ and, if $v_{\alpha}(x)<+\infty,$ then
$A_{\alpha}(x)$ is compact.
\item[{\rm(v)}] for an infinite-horizon $T=\infty$ there exists a stationary
discount-optimal policy $\phi_\alpha,$ and a stationary policy is
optimal if and only if $\phi_\alpha(x)\in A_\alpha(x)$ for all
$x\in \X.$
%
\end{itemize}
\end{theorem}
\begin{proof}
See Appendix \ref{app:B}.
\end{proof}

\begin{remark}\label{rem:new}{\rm Let us consider an MDP $(\X,\A,q,c)$ with
$\X=\X_W\times \X_Y,$ the stochastic kernel $q$ on $\X$ given $\X\times\A$ being semi-uniform Feller, and the nonnegative function $c:\X\times \A{\to} \overline{\mathbb{R}}$ being  $\K$-inf-compact. Then, Lemma~\ref{cor:Corollary 5.15.} implies that the stochastic kernel $q$ on $\X$ given $\X\times\A$ is weakly continuous. Therefore, \cite[Theorem 2]{FKNMOR} implies all assumptions and conclusions of Theorem~\ref{prop:dcoe} and, in addition, the functions $v_{t,\alpha}(\cdot)$ 
and $v_\alpha(\cdot)$ are lower semi-continuous for all $t=0,1,\ldots$ and $\alpha\ge 0$.}
\end{remark}

 We also remark that, if the cost function $c$ is nonnegative, then  optimality equations hold and stationary (Markov) optimal policies satisfy them for problems with an infinite (finite) horizons without any continuity assumptions on the transition probabilities $q$ and cost function $c;$  see, e.g.,  \cite[Propositions 9.8, 9.12 and Corollary 9.12.1]{BS} for $\alpha=1.$  This is also true, in the following two cases: (a) $c, \alpha\ge 0,$ and (b) $c\ge K>-\infty$ and $\alpha\in [0,1).$  However, if transition probabilities and costs do not satisfy appropriate continuity assumptions, then $\min$ should be replaced with $\inf$   in the optimality equations stated in statements (ii) and (iv) of Theorem~\ref{prop:dcoe}, the sets $A_{t,\alpha}(x)$ and $A_\alpha(x)$ can be empty, optimal policies may not exist, and, though a limit of value iterations with zero terminal costs exists, it may not be equal to the value function; see Yu~\cite{Yuj} and references therein on value iterations for infinite-state MDPs.

\section{Total-Cost Optimal Policies for MDPII and Corollaries for Platzman's  Model and for POMDPs}\label{sec:POMDPappl}

In this section we formulate Theorems~\ref{th:mainMDPII} and \ref{th:mainMDPIINew} stating the equivalences of semi-uniform Feller continuities of the transition probability $P$ for an MDPII,  stochastic kernel $R$ defined in \eqref{3.3}, and  transition probability $q$ for the MDPCI defined in \eqref{3.7}.  These two theorems also provide other necessary and sufficient conditions for semi-uniform Feller continuity of the stochastic kernels $P,$ $R$, and $q.$   The proofs of Theorems~\ref{th:mainMDPII} and \ref{th:mainMDPIINew} use Theorems~\ref{th:extra}, \ref{th:CBRmain}, the reduction of MDPIIs to MDPCIs established in \cite{Rh, Yu} and described in Section~\ref{subsec:Reduction}, and \cite[Theorem 3.3]{FKZ} stating that integration of cost functions with respect to probability measures in the argument corresponding to unobservable state variables preserves $\K$-inf-compactness of cost functions. Then we consider Platzman's  model and POMDPs and describe sufficient conditions for weak continuity of transition kernels in the reduced models, whose states are belief probabilities, and the validity of optimality equations, convergence of value iterations, and existence of optimal policies for these models.

\begin{theorem}\label{th:mainMDPII}
Let $(\W\times\Y,\A,P,c)$ be an MDPII, $(\P(\W)\times\Y,\A,q,\c)$ be its MDPCI, and $y\in\Y.$ Then the following conditions are equivalent:
\begin{itemize}
\item[{\rm(a)}] Assumption~\ref{AssKern} holds with $\S_1:=\W,$ $\S_2:=\Y,$ $\S_3:=\W\times\A,$ and $\Psi:=P(\,\cdot\,|\,\cdot\,,y,\,\cdot\,);$
\item[{\rm(b)}] the stochastic kernel $P(\,\cdot\,|\,\cdot\,,y,\,\cdot\,)$ on $\W\times\Y$ given
$\W\times\A$ is semi-uniform Feller;
\item[{\rm(c)}] the stochastic kernel $R(\,\cdot\,|\,\cdot\,,y,\,\cdot\,)$ on $\W\times\Y$ given $\P(\W)\times \A$ is semi-uniform Feller;
\item[{\rm(d)}]the marginal kernel $R(\W,\,\cdot\,|\,\cdot\,,y,\,\cdot\,)$ on $\Y$ given ${\P(\W)\times\A}$ is continuous in total variation, and the stochastic kernel $H(\,\cdot\,|\,\cdot\,,y,\,\cdot\,,\,\cdot\,)$ on $\W$ given
$\P(\W)\times \A\times\Y$ defined in \eqref{3.4} satisfies Assumption~\ref{Ass:H};
\item[{\rm({e})}] the stochastic kernel $q(\,\cdot\,|\,\cdot\,,y,\,\cdot\,)$ on $\P(\W)\times\Y$
given $\P(\W)\times \A$ is semi-uniform Feller.
\end{itemize}
Moreover, if nonnegative function $c$ is {$\M\K(\W\times \A,\Y)$}-inf-compact, and
for each $y\in\Y$ anyone of the above conditions~(a)--(e) holds,
then all the assumptions and  conclusions of Theorem~\ref{prop:dcoe} hold for the MDPCI $(\P(\W)\times\Y,\A,q,\c).$
\end{theorem}

\begin{theorem}\label{th:mainMDPIINew}
Let $(\W\times\Y,\A,P,c)$ be an MDPII, and $(\P(\W)\times\Y,\A,q,\c)$ be its MDPCI. Then the following conditions are equivalent:
\begin{itemize}
\item[{\rm(a)}] Assumption~\ref{AssKern} holds with $\S_1:=\W,$ $\S_2:=\Y,$ $\S_3:=\W\times\Y\times\A,$ and $\Psi:=P;$
\item[{\rm(b)}] the stochastic kernel $P$ on $\W\times\Y$ given
$\W\times\Y\times\A$ is semi-uniform Feller;
\item[{\rm(c)}] the stochastic kernel {$R$ on $\W\times\Y$ given $\P(\W)\times\Y\times \A$ is  semi-uniform Feller;}
\item[{\rm(d)}]{the marginal kernel $R(\W,\,\cdot\,|\,\cdot\,)$ on $\Y$ given $\P(\W)\times\Y\times\A$ is continuous in total variation, and} the stochastic kernel $H$ on $\W$ given
$\P(\W)\times\Y\times \A\times\Y$ defined in \eqref{3.4} satisfies Assumption~\ref{Ass:H};
\item[{\rm({e})}] the stochastic kernel $q$ on $\P(\W)\times\Y$
given $\P(\W)\times\Y\times \A$ is semi-uniform Feller.
\end{itemize}
Moreover, if the nonnegative function $c$ is $\K$-inf-compact, and
anyone of the above conditions~(a)--(e) holds, then all the assumptions and  conclusions of Theorem~\ref{prop:dcoe} hold for the MDPCI $(\P(\W)\times\Y,\A,q,\c)$, and the functions $v_{t,\alpha},$ $t=0,1,\ldots,$ and $v_\alpha$ are lower semi-continuous  on $\X$.
\end{theorem}

The proofs of Theorems~\ref{th:mainMDPII} and~\ref{th:mainMDPIINew} are provided in Appendix \ref{app:B}. We recall that $c,\alpha\ge 0$ in Theorems~\ref{prop:dcoe} and \ref{th:mainMDPII}. If $0\le\alpha<1$  and the function $c$ is bounded below, then all conclusions of Theorems~\ref{prop:dcoe} and \ref{th:mainMDPII} hold with the following minor modifications (i)  the functions $v_{t,\alpha}$ and $v_\alpha$ are bounded below rather than nonnegative, and (ii) $v_{t,\alpha}(x)\to
v_\alpha (x)$ rather than $v_{t,\alpha}(x)\uparrow
v_\alpha (x)$ as $t\to\infty.$  This is true for function $c$ bounded below by $-K>-\infty$ because such MDPII can be converted into a model with nonnegative costs by replacing costs $c$ with  $c+K;$ \cite{FKZ}. The suggestion to fix $y$ in assumptions of Theorems~\ref{prop:dcoe} and \ref{th:mainMDPII} was proposed by a referee.

%

According to \cite{Rh, Yu}, for each optimal policy for the MDPCI $(\P(\W)\times\Y,\A,q,\c)$ there constructively exists an optimal policy in the original MDPII $(\W\times\Y,\A,P,c)$.
\cite[Theorem~4.4]{Steklov} establishes weak continuity of the transition kernel in the MDPCI under the more restrictive assumption than statement~(a) of Theorem~\ref{th:mainMDPII} when the countable base in {Assumption~\ref{AssKern}} does not depend on the argument $s_3=(w,y,a);$ see also \cite[Example~1]{FKZJThProb}. Moreover, for any $T=1,2,\ldots$ and $\alpha\ge0,$ the value functions $\tilde{V}_{T,\alpha}(z,y),\tilde{V}_\alpha(z,y)$ in the MDPCI $(\P(\W)\times\Y,\A,q,\c)$  are concave in $z\in \P(\W).$  This is true because infimums of affine functions are concave functions.


The proof of Theorem~\ref{th:mainMDPII} uses the following preservation property for {$\M\K(\W\times\A,\Y)$}-inf-compactness.

\begin{theorem}\label{th:PresC}
If $c:\,\W\times\Y\times\A{\to}  \overline{\R}_+$ is an {$\M\K(\W\times\A,\Y)$}-inf-compact function, then the function $\c:\P(\W)\times\Y\times\A{\to}\overline{\R}_{+}$ defined in \eqref{eq:c} is {$\M\K(\P(\W)\times\A,\Y)$}-inf-compact.
\end{theorem}

\begin{proof} This theorem
follows from \cite[Proposition~7.29]{BS} on preservation of Borel measurability and from \cite[Theorem~3.3]{FKZ} on preservation of $\K$-inf-compactness.
\end{proof}

The particular case of an MDPII is a probabilistic dynamical system considered in Platzman \cite{Plat}.
\begin{definition}\label{defi:Plpl}
Platzman's  model is specified by an MDPII $(\W\times\Y,\A,P,c),$ where $P$ is a stochastic kernel on $\W\times\Y$ given
$\W\times\A.$
\end{definition}
\begin{remark}\label{rem:Plpl}
{\rm
Formally speaking, Platzman's  model is an MDPII with the transition kernel $P(\,\cdot\,|w,y,$ $a)$ that does not depend
on $y$. Therefore, Theorem~\ref{th:mainMDPII} implies certain corollaries for Platzman's  model.
}
\end{remark}

\begin{corollary}\label{cor:Plpl}
Let $(\W\times\Y,\A,P,c)$ be Platzman's  model. Then the stochastic kernel $P$ on $\W\times\Y$ given
$\W\times\A$ is semi-uniform Feller if and only if one of the equivalent conditions~(a), (c), {(d), or (e)} of Theorem~\ref{th:mainMDPII} holds.
Moreover, if the nonnegative function $c$ is {$\M\K(\W\times\A,\Y)$-inf-compact and the stochastic kernel} $P$ on $\W\times\Y$ given $\W\times\A$ is semi-uniform Feller, then all {the assumptions and conclusions} of Theorem~\ref{th:mainMDPII} hold.
\end{corollary}
\begin{proof}
According to Remark~\ref{rem:Plpl}, Corollary~\ref{cor:Plpl} follows directly from Theorem~\ref{th:mainMDPII}.
\end{proof}

For Platzman's models we shall write $P(B\times C|w,a),$ $R(B\times C|z,a),$ $H(D|z,a,y'),$ and $q(D\times C|z,a)$ instead of  $P(B\times C|w,y,a),$ $R(B\times C|z,y,a),$ $H(D|z,y{,}a,y'),$ and $q(D\times C|z,y,a)$ since these stochastic kernels do not depend on the variable  $y.$ For Platzman's models we shall also consider the marginal kernel
$\hat{q}(D|z,a):={q(D,\Y|z,a)}$ on $\P(\W)$ given $\P(\W)\times\A.$  In view of \eqref{3.7}, for $(z,a)\in \P(\W)\times\A$
and for $D\in \B(\P(\W)),$
\begin{equation}\label{3.777A}
\begin{aligned} \hat{q}(D|z,a):=\int_{\Y}\h\{H(z,a,y')\in D\}
R(\W,dy'|z, a).\end{aligned}
\end{equation}

\begin{corollary}\label{cor:Plpl1}
Let $(\W\times\Y,\A,P,c)$ be Platzman's  model, and let the stochastic kernel $P$ on $\W\times\Y$ given
$\W\times\A$ {be} semi-uniform Feller. Then the stochastic kernel $\hat{q}$ on $\P(\W)$ given
$\P(\W)\times\A$ is weakly continuous.
\end{corollary}

\begin{proof}
According to Corollary~\ref{cor:Plpl} {and Lemma~\ref{cor:Corollary 5.15.},} the stochastic kernel $q$ on $\P(\W)\times \Y$ given   $\P(\W)\times \A$ is weakly continuous.  Therefore, its marginal kernel $\hat{q}$ on $\P(\W)$ given   $\P(\W)\times \A$ is also weakly continuous.
\end{proof}

As mentioned in {\cite{Plat}, the special cases of Platzman's}  model include  two partially observable MDPs  which we denote as ${\rm POMDP}_1$ and ${\rm POMDP}_2;$ see Definitions~\ref{defi:POMDP1}, \ref{defi:POMDP2} and Figure~\ref{Fig:Diagram}.
%

Let $i=1,2,$ let $\W,$ $\Y,$ and $\A$ be Borel subsets of Polish spaces, $P_i(dw'|w,a)$ be a stochastic kernel on
$\W$ given $\W\times\A,$ $Q_1(dy| w,a)$ be a stochastic kernel on
$\Y$ given $\W\times\A,$ $Q_2(dy| a,w)$ be a stochastic kernel on
$\A$ given $\A\times\W,$ $Q_{0,i}(dy|w)$ be a stochastic kernel on
$\Y$ given $\W,$ $p$ be a probability distribution on $\W.$
\begin{definition}\label{defi:POMDP1} A
${\rm POMDP}_1$ $(\W,\Y,\A,P_1,Q_1,c)$ is specified by Platzman's model $(\W\times\Y,\A,P,c)$ with \begin{equation}\label{eq:neq1}
P(B\times C| w,a):=P_1(B|w,a)Q_1(C|w,a),
\end{equation}
$B\in\B(\W),$ $C\in\B(\Y),$ $w\in\W,$ $y\in\Y,$ $a\in\A.$
\end{definition}

Let $(\W,\Y,\A,P_1,Q_1,c)$ be a ${\rm POMDP}_1.$ Then, the stochastic kernel $R$ on $\W\times\Y$ given
$\P(\W)\times \A,$  which is  defined for MDPIIs in \eqref{3.3}, takes the following form,
\begin{equation}\label{eq:3.3POMDP1}
R(B\times C|z,a):=\int_{\W} Q_1(C|w,a)P_1(B|w,a)z(dw),
\end{equation}
$B\in \mathcal{B}(\W),$ $C\in \mathcal{B}(\Y),$ $z\in\P(\W),$ $a\in \A.$
\begin{definition}\label{defi:POMDP2} A
${\rm POMDP}_2$ $(\W,\Y,\A,P_2,Q_2,c)$ is specified by Platzman's model $(\W\times\Y,\A,P,c)$ with \begin{equation}\label{eq:defi2}
P(B\times C| w,a):=\int_B Q_2(C|a,w')P_2(dw'|w,a),
\end{equation}
$B\in\B(\W),$ $C\in\B(\Y),$ $w\in\W,$ $y\in\Y,$ $a\in\A.$
\end{definition}

 We recall that Figure~\ref{Fig:Diagram} describes the relations between an MDPII, Platzman's model, ${\rm POMDP}_1,$  and ${\rm POMDP}_2$  based  on the generality of transition probabilities $P.$ In addition, ${\rm POMDP}_1$ and ${\rm POMDP}_2$ are two different models.  For example, for a ${\rm POMDP}_1$ the random variables $w_{t+1}$ and $y_{t+1}$ are conditionally independent  given the values $w_t$ and $a_t.$  This is not true for  ${\rm POMDP}_2.$

 Other relations between these models  also take place.  In particular, a reduction of {an MDPII to a ${\rm POMDP}_2$   is described in \cite[Section 6]{Steklov} and in \cite[Section 8.3]{FKZ}.   Therefore, in some sense
an MDPII,\ Platzman's model, and  a ${\rm POMDP}_2$ can be viewed as equivalent models. This reduction was used in \cite{FKZ}} to    prove Theorem 8.1 there stating sufficient
conditions for weak continuity of transition probabilities for MDPCIs. {This reduction transforms an MDPII with a weakly continuous transition
probability  into a  ${\rm POMDP}_2$ with weakly continuous transition and observation probabilities.  Since weak continuity of transition and observation probabilities for ${\rm POMDP}_2$ are not sufficient for continuity
of transition probabilities for the corresponding belief-MDP (see \cite[Example 4.1]{FKZ}), \cite[Theorem 8.1]{FKZ} contains an additional assumption on the transition probability $P$ of the MDPII.
This assumption is relaxed in \cite[Theorem 6.2] {Steklov}.}
As shown in \cite[Example 1]{FKZJThProb},   semi-uniform Feller continuity of the transition probability $P$ assumed in this paper is a more general property than the assumption on $P$ in \cite[Theorem 6.2] {Steklov}.

For a ${\rm POMDP}_2$ $(\W,\Y,\A,P_2,Q_2,c)$ the stochastic kernel $R$ on $\W\times\Y$ given $\P(\W)\times \A,$  which is defined for MDPIIs in \eqref{3.3}, takes the following form,
\begin{equation}\label{eq:3.3POMDP2}
R(B\times C|z,a):=\int_{\W}\int_B Q_2(C|a,w')P_2(dw'|w,a)z(dw),
\end{equation}
$B\in \mathcal{B}(\W),$ $C\in \mathcal{B}(\Y),$ $z\in\P(\W),$ $a\in \A.$
A ${\rm POMDP}_1$ is Platzman's model with observations $y_{t+1}$ being ``random functions'' of $w_{t}$ and $a_{t},$ and a ${\rm POMDP}_2$ is Platzman's  model with observations
$y_{t+1}$ being ``random functions'' of $a_{t}$ and $w_{t+1}.$
Let us apply Theorem~\ref{th:mainMDPII} to a ${\rm POMDP}_1$ and ${\rm POMDP}_2.$


{\rm
Corollary~\ref{cor:Plpl} establishes necessary and sufficient conditions for semi-uniform Feller continuity of the transition probabilities $P$ for Platzman's model $(\W\times\Y,\A,P,c)$ in terms of the same property for the transition probabilities $q$ of the respective belief-MDP $(\P(\W)\times\Y,\A,q,\c).$ Since a ${\rm POMDP}_i,$ $i=1,2,$ is a particular case of Platzman's model, Corollary~\ref{cor:Plpl} implies the necessary and sufficient conditions for semi-uniform Feller continuity of the stochastic kernel $q$ on $\Y\times\P(\W)$ given $\W\times\A$ in terms of the same property for the transition probability $P$ defined in \eqref{eq:neq1} for a ${\rm POMDP}_1$ and in \eqref{eq:defi2} for a ${\rm POMDP}_2$ respectively.
}
%
\begin{corollary}\label{cor:POMDPmain1}
For a ${\rm POMDP}_1$ $(\W,\Y,\A,P_1,Q_1,c)$, the following two conditions holding together:
\begin{itemize}\item[{\rm(a)}] the stochastic kernel $P_1$ on
$\W$ given $\W\times\A$ is weakly continuous; \item[{\rm(b)}] the stochastic kernel $Q_1$ on
$\Y$ given $\W\times\A$ is continuous in total variation;
 \end{itemize} are equivalent to semi-uniform Feller continuity of the stochastic kernel $P$ on $\W\times\Y$ given $\W\times\A.$ Moreover, if these two conditions hold, then:
\begin{itemize}\item[{\rm(i)}]
statements~{(a), (c)}--({e}) of Theorem~\ref{th:mainMDPII} hold;
 \item[{\rm(ii)}] if the nonnegative function $c:\W\times\Y\times\A{\to}\overline{\R}$ is {$\M\K(\W\times\A,\Y)$}-inf-compact, then all the conclusions of Theorem~\ref{th:mainMDPII} hold;
     \item[{\rm(iii)}]the stochastic kernel $\hat{q}$ on $\P(\W)$ given
$\P(\W)\times\A$ defined in \eqref{3.777A} is weakly continuous.\end{itemize}
\end{corollary}
\begin{proof}
See Appendix~\ref{app:B}.
\end{proof}
\begin{corollary}\label{cor:POMDPmain}
For a ${\rm POMDP}_2$ $(\W,\Y,\A,P_2,Q_2,c)$ each of the following conditions:
\begin{itemize}
\item[{\rm(a)}] the stochastic kernel $P_2$ on $\W$ given $\W\times\A$ is
weakly continuous, and the stochastic kernel $Q_2$ on $\Y$
given $\A\times\W$ is continuous in
 total  variation;
\item[{\rm(b)}] the stochastic kernel $P_2$ on $\W$ given $\W\times\A$ is
continuous in total  variation, and the observation kernel $Q_2$ on $\Y$ given $\A\times\W$ is continuous in $a$ in total  variation;
\end{itemize}
implies semi-uniform Feller continuity of the stochastic kernel $P$ on $\W\times\Y$ given $\W\times\A.$ 
Moreover, each of conditions (a) or (b) implies the validity of conclusions (i)--(iii) of Corollary~\ref{cor:POMDPmain1}  for the  ${\rm POMDP}_2$.  
\end{corollary}
\begin{proof}
See Appendix~\ref{app:B}.
\end{proof}

Regarding Corollary~\ref{cor:POMDPmain}, weak continuity of the stochastic kernel $\hat{q}$ on $\P(\W)\times \A$ for a ${\rm POMDP}_2$ under condition (a) from Corollary~\ref{cor:POMDPmain} is stated in \cite[Theorem 3.6]{FKZ}, and another proof of this statement is  provided in \cite[Theorem 1]{Saldi}. Weak continuity of the stochastic kernel $\hat{q}$ on $\P(\W)\times \A$ for a ${\rm POMDP}_2$ under condition (b) from Corollary~\ref{cor:POMDPmain} is an extension of \cite[Theorem 2]{Saldi}, where this weak continuity is proved under the assumption that the stochastic kernel $P_2$ on $\W$ given $\W\times\A$ is continuous in total variation and the observation kernel $Q_2$ does not depend on actions.

{Different sufficient conditions for weak continuity of the kernel  $\hat{q}$ for a ${\rm POMDP}_2$ are formulated in  monographs \cite{HL} and  \cite{RS}.   In both cases these conditions are stronger than condition (a) from Corollary~\ref{cor:POMDPmain}.  In terms of the current paper,   weak continuity of the stochastic kernel $\hat{q}$ on $\P(\W)$ given $\P(\W)\times\A$ is stated in \cite[p. 92]{HL} under condition (a) from Corollary~\ref{cor:POMDPmain} and under the assumption that the observation space $\Y$ is denumerable.  The proof on \cite[p. 93]{HL} is based on the existence of a transition kernel $H(z,a,y'),$ which is weakly continuous in $(z,a,y')$ and satisfies \eqref{3.777A}.  However, \cite[Example 4]{FKZ} shows that such kernel may not exist { even} for a ${\rm POMDP}_2$ with finite sets $\X,$ $\Y$ and {continuous in $a$} functions $P_2(x'|x,a)$ and $Q_2(y|a,x).$  A ${\rm POMDP}_2$ is considered in \cite[Chapter 2]{RS} under additional assumptions that the state space $\X$ is locally compact, observations $y_t$ belong to an Euclidean space, and the observation kernel does not depend on actions and has a density, that is, $Q(dy|x)=r(x,y)dy.$  Weak continuity of the kernel  $\hat{q}$ is stated in \cite[Corollary 1.5]{RS} under four assumptions, which taken together are stronger than condition (a) in Corollary~\ref{cor:POMDPmain}.}

Let us consider Platzman's model $(\W\times\Y,\A,P,c)$ with the cost function $c$ that does not depend on observations $y,$ that is, $c(w,y,a)=c(w,a).$  In this case the MDPCI $(\P(\W)\times\Y,\A,q,\c)$ can be reduced to a smaller MDP $(\P(\W),\A,{\hat q},{\hat c})$ with the state space $\P(\W),$ action space $\A,$ transition probability ${\hat q}$ defined in \eqref{3.777A}, and
 one-step cost  function ${\hat c}:\P(\W)\times\A{\to}\overline{\R}$, defined for $z\in\P(\W)$ and $a\in\A$ as
\begin{equation}\label{eq:chat}
\hat{c}(z,a):=\int_{\W}c(w,a)z(dw).
\end{equation}
The reduction of an MDPCI $(\P(\W)\times\Y,\A,q,\c)$ to the belief-MDP $(\P(\W),\A,{\hat q},{\hat c})$ holds in view of \cite[Theorem~2]{F2005} because in the MDPCI transition probabilities from  states $(z_t,y_t)\in \P(\W)\times\Y$ to states $z_{t+1}\in\P(\W)$  and  costs $c(z_t,a_t)$ do not depend on $y_t.$  If a Markov or stationary optimal policy is found for the belief-MDP $(\P(\W),\A,{\hat q},{\hat c}),$ it is possible, { as described at the end of Section~\ref{subsec:Reduction},} to construct an optimal policy for Platzman's models following the same procedures as constructing an optimal policy for and MDPII given a Markov or stationary optimal policy for the corresponding MDPCI.

\begin{corollary}\label{cor:POMDPmaincwa}Let us consider Platzman's model $(\W\times\Y,\A,P,c)$ with the one-step cost function $c:\W\times\A\to\overline{\R}_+ .$ If  the stochastic kernel $P$ on $\W\times\Y$ given
$\W\times\A$ is semi-uniform Feller, and the one-step cost function $c$ is $\K$-inf-compact on $\W\times\A$, then the transition kernel $\hat{q}$ on $\P(\W)$ given $\P(\W)\times \A$ is weakly continuous, the one-step cost function $\hat{c}$ is $\K$-inf-compact on $\P(\W)\times\A$, and all the conclusions of \cite[Theorem 2.1]{FKZ} hold for the belief-MDP  $(\P(\W),\A,{\hat q},{\hat c}),$ that is:
 \begin{itemize} \item[{\rm(i)}]
 optimality equations hold, {and they define optimal policies;}
 \item[{\rm(ii)}] value iterations converge {to optimal values} if zero terminal costs are chosen;
 \item[{\rm(iii)}]    Markov optimal policies exist for finite-horizon problems;
 \item[{\rm(iv)}]    stationary optimal policies exist for infinite-horizon problems.
 \end{itemize}
 {Moreover, all {these conclusions} hold for {a} ${\rm POMDP}_1$ $(\W,\Y,\A,P_1,Q_1,c)$ with the transition and observation kernels $P_1$ and $Q_1$ satisfying conditions (a) and (b) from Corollary~\ref{cor:POMDPmain1} and for {a} ${\rm POMDP}_2$ $(\W,\Y,\A,P_2,Q_2,c)$ with the transition and observation  kernels $P_2$ and $Q_2$ satisfying either condition (a) or condition (b)  from Corollary~\ref{cor:POMDPmain}.}
\end{corollary}
\begin{proof}
Weak continuity of the stochastic kernel $\hat{q}$  on $\P(\W)$ given $\P(\W)\times \A$ is stated in Corollary~\ref{cor:Plpl1}. $\K$-inf-compactness of the function $\hat{c}$ on $\P(\W)\times\A$ follows from \cite[Theorem 3.3]{FKZ}.  The remaining statements of the corollary follow from    \cite[Theorem 2.1]{FKZ}. {The transition probability $P$ for ${\rm POMDP}_1$ $(\W,\Y,\A,P_1,Q_1,c)$ defined in \eqref{eq:neq1} is semi-uniform Feller according to Corollary~\ref{cor:POMDPmain1}, and the transition probability $P$ for ${\rm POMDP}_2$ $(\W,\Y,\A,P_2,Q_2,c)$ defined in \eqref{eq:defi2} is semi-uniform Feller due to Corollary~\ref{cor:POMDPmain}.}
\end{proof}

\appendix

\section{Proofs of Theorems~\ref{th:CBRmain}, \ref{prop:dcoe}, \ref{th:mainMDPII}, and Corollaries~\ref{cor:POMDPmain1}, \ref{cor:POMDPmain}}\label{app:B}

{
We use the following fact in the proofs of equalities \eqref{eq:neww01} and \eqref{eq:neww1} below: if $\{G^{\left(n\right)},G\}_{n=1,2,\ldots}$ {is} a sequence of finite measures on a metric space $\mathcal{S}$ and $\{g^{\left(n\right)},g\}_{n=1,2,\ldots}$ {is} a uniformly bounded sequence of Borel measurable functions on $\mathcal{S}$ such that
\[
\lim_{n\to\infty}\sup_{B\in\B(\mathcal{S})}\left|\int_{B}g^{\left(n\right)}(s)G^{\left(n\right)}(ds) - \int_{B}g^{\left(n\right)}(s)G(ds) \right|=0,
\]
then
\[
\lim_{n\to\infty}\sup_{B\in\B(\mathcal{S})}\left|\int_{B}g^{\left(n\right)}(s)G^{\left(n\right)}(ds) - \int_{B}g(s)G(ds) \right|=0
\]
holds if and only if
\[
\lim_{n\to\infty}\sup_{B\in\B(\mathcal{S})}\left|\int_{B}g^{\left(n\right)}(s)G(ds) - \int_{B}g(s)G(ds) \right|=0.
\]
}

\begin{proof}[Proof of Theorem~\ref{th:CBRmain}]
(a) $\Rightarrow$ (b). {Since the stochastic kernel $\Psi$ on $\S_1\times\S_2$ given $\S_3$ is semi-uniform Feller, the marginal kernel $\Psi(\S_1, \,\cdot\,|\,\cdot\,)$ is continuous in total variation. Moreover, for each bounded continuous function $f$ on $\S_1,$ we have from  \eqref{eq:equivWTV3} and \eqref{eq:CBR1} that
\begin{equation}\label{eq:neww01}
\lim_{n\to\infty} \sup_{B\in \B(\S_2)} \left| \int_{B}\int_{\S_1} f(s_1)\Phi(ds_1|s_2,s_3^{\left(n\right)})\Psi(\S_1,ds_2|s_3)-\int_{B}\int_{\S_1} f(s_1)\Phi(ds_1|s_2,s_3)\Psi(\S_1,ds_2|s_3)\right|= 0
\end{equation}
because the family of Borel measurable functions $\{s_2\mapsto\int_{\S_1} f(s_1)\Phi(ds_1|s_2,s_3^{\left(n\right)})\,:\,n=1,2,\ldots\}$ is uniformly bounded on $\S_2$ by the same constant as $f$ on $\S_1.$ This is equivalent to $\int_{\S_1} f(s_1)\Phi(ds_1|\,\cdot\,,s_3^{\left(n\right)})\to \int_{\S_1} f(s_1)\Phi(ds_1|\,\cdot\,,s_3)$ in $L_1(\S_2,\B(\S_2),\nu)$ with $\nu(\,\cdot\,):= \Psi(\S_1,\,\cdot\,|s_3).$ {Therefore,}
\[
\int_{\S_1} f(s_1)\Phi(ds_1|\,\cdot\,,s_3^{\left(n_k\right)})\to \int_{\S_1} f(s_1)\Phi(ds_1|\,\cdot\,,s_3) \quad \nu\mbox{-almost surely, as }k\to \infty,
\]
for some sequence $\{n_k\}_{k=1,2,\ldots}$ ($n_k\uparrow\infty$ as $k\to\infty$). { We apply the diagonalization procedure} to extract a subsequence $\{\tilde{n}_k\}_{k=1,2,\ldots}$ ($\tilde{n}_k\uparrow\infty$ as $k\to\infty$) such that
\[
\int_{\S_1} g(s_1)\Phi(ds_1|\,\cdot\,,s_3^{\left(\tilde{n}_k\right)})\to \int_{\S_1} g(s_1)\Phi(ds_1|\,\cdot\,,s_3) \quad \nu\mbox{-almost surely, as }k\to\infty,
\]
for each $g\in \mathcal{G},$ where $\mathcal{G}$ is a countable uniformly bounded  family of continuous functions on $\S_2$ that determines weak convergence of probability measures on $\S_2$ according to Parthasarathy \cite[Theorem~6.6, p.~47]{Part}. Thus, $\Phi(\,\cdot\,,s_3^{\left(\tilde{n}_k\right)})$ converges weakly to $\Phi(\,\cdot\,,s_3)$ $\nu$-almost surely, {and}  Assumption~\ref{Ass:H} holds.
}

({b}) $\Rightarrow$ ({c}).
Let $f$ be a bounded continuous function on $\P(\S_1).$ Since $\Psi(\S_1, \,\cdot\,|\,\cdot\,)$ is continuous in total variation, to prove that \eqref{eq:equivWTV3} holds for the stochastic kernel $\phi,$ it is sufficient to show that
\begin{equation}\label{eq:neww1}
\lim_{n\to\infty} \sup_{B\in \B(\S_2)} \left| \int_{B} f(\Phi(s_2,s_3^{\left(n\right)})) \Psi(\S_1,ds_2|s_3)-\int_{B} f(\Phi(s_2,s_3)) \Psi(\S_1,ds_2|s_3)\right|= 0.
\end{equation}
For the probability space $\Sigma:=(\S_2,\B(\S_2),\mu)$ with $\mu(\,\cdot\,):=\Psi(\S_1,\,\cdot\,|s_3),$ the $\P(\S_1)$-valued random variables $\Phi(\,\cdot\,,s_3^{\left(n\right)})\mathop{\to}\limits^\mu\Phi(\,\cdot\,,s_3)$ as $n\to\infty,$ according to Assumption~\ref{Ass:H} and \eqref{eq:CBR4}, where $\nu^{(n)}\mathop{\to}\limits^\mu \nu$ denotes the convergence in probability $\mu,$ that is,  $\rho_{\P(\S_1)}(\nu^{(n)},\nu)\to 0$ in probability $\mu.$ Then $f(\Phi(\,\cdot\,,s_3^{\left(n\right)}))\mathop{\to}\limits^\mu f(\Phi(\,\cdot\,,s_3))$ because $f$ is continuous on $\P(\S_1).$ In turn, since $f$ is bounded on $\P(\S_1),$ this implies that $f(\Phi(\,\cdot\,,s_3^{\left(n\right)}))\to f(\Phi(\,\cdot\,,s_3))$ in $L_1(\Sigma),$ from which the desired relation \eqref{eq:neww1} follows.

({c}) $\Rightarrow$ (a). Let a sequence
$\{s_3^{\left(n\right)}\}_{n=1,2,\ldots}\subset\S_3$ converge to
$s_3\in\S_3$ as $n\to\infty.$ {Since the stochastic kernel $\phi$ on $\P(\S_1)\times \S_2$ given $\S_3$ is semi-uniform Feller, for every nonnegative bounded lower semi-continuous function $f$ on $\P(\S_1)$, according to Theorem~\ref{th:equivWTV}(a,e),
\begin{equation}\label{eqEF1129}
\ilim_{n\to\infty} \inf_{B\in \B(\S_2)} \left( \int_{\P(\S_1)} f(\mu) \phi(d\mu,B|s_3^{\left(n\right)})-\int_{\P(\S_1)} f(\mu) \phi(d\mu,B|s_3)\right)= 0.
\end{equation}

 For each $B\in \B(\S_2),$ formula \eqref{eq:CBR2} establishes the equality of two measures
on $(\P(\S_1),\B(\P(\S_1))).$
Therefore,
 for every Borel measurable nonnegative functions $f$ on $\P(\S_1),$
\begin{equation}\label{eqEF1129A}
 \int_B f(\Phi(s_2,\tilde{s}_3)) \Psi(\S_1,ds_2|\tilde{s}_3)=\int_{\P(\S_1)}f(\mu) \phi(d\mu,B|\tilde{s}_3),\qquad \tilde{s}_3\in\S_3.
\end{equation}

Let us fix an arbitrary
open set $\oo\subset \S_1$ and consider nonnegative bounded lower semi-continuous function $f(\mu):=\mu(\oo),$ $\mu\in\P(\S_1).$ Then
\[
\begin{aligned}
\ilim_{n\to\infty}&\inf_{B\in\B(\S_2)}\left(\int_{B} \Phi(\oo|s_2,s_3^{\left(n\right)})\Psi(\S_1,ds_2|s_3^{\left(n\right)})
 -\int_{B}\Phi(\oo|s_2,s_3)\Psi(\S_1,ds_2|s_3)
\right)\\
&=\ilim_{n\to\infty} \inf_{B\in \B(\S_2)} \left( \int_B f(\Phi(s_2,s_3^{\left(n\right)})) \Psi(\S_1,ds_2|s_3^{\left(n\right)})-\int_B f(\Phi(s_2,s_3)) \Psi(\S_1,ds_2|s_3)\right)=0,
\end{aligned}
\]
where the first equality follows from the definition of $f$, and the second equality follows from \eqref{eqEF1129A} and from \eqref{eqEF1129}. Thus, the stochastic kernel $\Psi$ on $\S_1\times\S_2$ given $\S_3$ is WTV-continuous, and therefore it is semi-uniform Feller.}
\end{proof}

\begin{remark}\label{rem:J}
{\rm
Theorem~\ref{th:CBRmain} can be proved in multiple ways using  equivalent characterizations of semi-uniform Feller kernels. The original proofs \cite[Proof of Theorem~5.10, pp.~16--20]{FetalOld} were based on some of these
 characterizations, while the current proofs of (a) $\Rightarrow$ (b) $\Rightarrow$ (c) were suggested by a referee.}
\end{remark}

The following Lemma~\ref{lm:NC1} is useful for establishing continuity
properties of the value functions $v_{n,\alpha}(x)$ and
$v_\alpha(x)$ in $x\in\X$ {stated in} Theorem~\ref{prop:dcoe}.

\begin{lemma}\label{lm:NC1} Let the MDP $(\X,\A,q,c)$ {satisfy} the assumptions of Theorem~\ref{prop:dcoe}, and let $\alpha\ge0{.}$ Then the function $u^*(x):=\inf\limits_{a\in \A}\eta_u^\alpha(x,a),$ $x\in \mathbb{X},$ where the function $\eta_u^\alpha$ is defined in \eqref{e:defeta}, belongs to $LW(\mathbb{X}),$ and there exists a stationary policy $f:\X{\to} \A$
such that $u^*(x):=\eta_u^\alpha(x,f(x)),$ $x\in \mathbb{X}.$
Moreover, the sets
$A_*(x)=\left\{a\in \A\,:\,u^*(x)=\eta_u^\alpha(x,a)\right\},$ $x\in \X,$ which are nonempty, satisfy the following properties: (a) the graph ${\rm Gr}(A_*)=\{(x,a)\,:\, x\in\X,\, a\in A_*(x)\}$ is a Borel subset of $\X\times \A;$ (b) if $u^*(x)=+\infty,$
then $A_*(x)=\A,$ and, if $u^*(x)<+\infty,$ then $A_*(x)$ is
compact.
\end{lemma}

\begin{proof}
The function $(x,a)\mapsto \eta_{u}^\alpha(x,a)$ is \textit{nonnegative} because $c,$ $u,$ and $\alpha$ are nonnegative. Therefore, since $u$ is a Borel measurable function, and   $q$ is a stochastic kernel,
\cite[Proposition~7.29]{BS} implies that the function $(x,a)\mapsto \int_\X  u(\tilde{x})q(d\tilde{x}|x,a)$ is Borel measurable on $\X\times\A,$ which implies that the function $(x,a)\mapsto \eta_{u}^\alpha(x,a)$ is \textit{Borel measurable} on $\X\times\A$ because $c$ is Borel measurable.

Let us prove that the function $(w,a)\mapsto \int_\X  u(\tilde{x})q(d\tilde{x}|w,y,a)$ is l.s.c.\ on $\X_W\times\A$ for each $y\in\X_Y.$ On the contrary, if this function is not l.s.c., then there exist a sequence $\{(w^{\left(n\right)},a^{\left(n\right)})\}_{n=1,2,\ldots}\subset \X_W\times\A$ converging to some $(w,a)\in\X_W\times\A$ and a constant $\lambda$ such that for each $n=1,2,\ldots$
\begin{equation}\label{eq:Int1}
\int_{\X_W\times\X_Y}  u(\tilde{w},\tilde{y})q(d\tilde{w}\times d\tilde{y}|w^{\left(n\right)},y,a^{\left(n\right)})\le \lambda<\int_\X  u(\tilde{x})q(d\tilde{x}|w,y,a).
\end{equation}
According to Theorem~\ref{th:CBRmain}(a,b) applied to $\Psi:=q,$ $\S_1:=\X_W,$ $\S_2:=\X_Y,$ $\S_3:=\X_W\times\{y\}\times\A,$ there exists a stochastic kernel $\Phi$ on $\X_W$ given $\X_Y\times\X_W\times\{y\}\times\A$ such that \eqref{eq:CBR1} and Assumption~\ref{Ass:H} hold. In particular, \eqref{eq:Int1} implies that for each $n=1,2,\ldots$
\[
\int_{\X_Y}\left[
\int_{\X_W}u(\tilde{w},\tilde{y})\Phi(d\tilde{w}|\tilde{y},w^{\left(n\right)},y,a^{\left(n\right)})
\right]q(\X_W,d\tilde{y}|w^{\left(n\right)},y,a^{\left(n\right)})\le \lambda,
\]
and there exist a subsequence $\{(w^{\left(n_k\right)},a^{\left(n_k\right)})\}_{k=1,2,\ldots}\subset\{(w^{\left(n\right)},a^{\left(n\right)})\}_{n=1,2,\ldots}$ and a Borel set $Y\in\B(\X_Y)$ such that $q(\X_W\times Y|w,y,a)=1$ and $\Phi(\tilde{y},w^{\left(n\right)},y,a^{\left(n\right)})$ converges weakly to $\Phi(\tilde{y},w,y,a)$ in $\P(\X_W)$ as $k\to \infty,$ for all $\tilde{y}\in Y.$ Therefore, since the function $\tilde{w}\mapsto u(\tilde{w},\tilde{y})$ is nonnegative and l.s.c.\ for each $\tilde{y}\in Y,$ Fatou's lemma for weakly converging probabilities \cite[Theorem~1.1]{TVP} implies that for each $\tilde{y}\in Y$
\begin{equation}\label{eqlastEF}
\int_{\X_W}u(\tilde{w},\tilde{y})\Phi(d\tilde{w}|\tilde{y},w,y,a)\le \ilim_{k\to\infty}\int_{\X_W}u(\tilde{w},\tilde{y})\Phi(d\tilde{w}|\tilde{y},w^{\left(n_k\right)},y,a^{\left(n_k\right)}).
\end{equation}
For a fixed $N=1,2,\ldots,$ we set $\varphi_k^N(\tilde{y}):=\min\{\int_{\X_W}u(\tilde{w},\tilde{y})\Phi(d\tilde{w}|\tilde{y},w^{\left(n_k\right)},y,a^{\left(n_k\right)}),N\}$ and $\varphi^N(\tilde{y}):=\min\{\int_{\X_W}u(\tilde{w},\tilde{y})\Phi(d\tilde{w}|\tilde{y},w,y,a),N\},$ where $\tilde{y}\in Y,$ $k=1,2,\ldots.$ Note that
$\varphi^N(\tilde{y})\le \ilim_{k\to\infty}\varphi_k^N(\tilde{y}),$ $\tilde{y}\in Y,$ {in view of \eqref{eqlastEF}}. Therefore, uniform Fatou's lemma \cite[Corollary~2.3]{UFL} implies that for each $N=1,2,\ldots$
\[
\begin{aligned}
\int_{\X_Y}&\varphi^N(\tilde{y})q(\X_W,d\tilde{y}|w,y,a)\le \ilim_{k\to\infty}\int_{\X_Y}\varphi_k^N(\tilde{y})q(\X_W,d\tilde{y}|w^{\left(n_k\right)},y,a^{\left(n_k\right)})\\
&\le \ilim_{k\to\infty}\int_{\X_Y}\left[
\int_{\X_W}u(\tilde{w},\tilde{y})\Phi(d\tilde{w}|\tilde{y},w^{\left(n_k\right)},y,a^{\left(n_k\right)})
\right]q(\X_W,d\tilde{y}|w^{\left(n_k\right)},y,a^{\left(n_k\right)})\le \lambda.
\end{aligned}
\]
 Thus, the monotone convergence theorem implies 
\[
 \int_\X  u(\tilde{x})q(d\tilde{x}|w,y,a) = \lim_{N\to \infty}\int_{\X_Y}\varphi^N(\tilde{y})q(\X_W,d\tilde{y}|w,y,a)\le \lambda.
\]
This is a contradiction with \eqref{eq:Int1}. Therefore, the function $(w,a)\mapsto \int_\X  u(\tilde{x})q(d\tilde{x}|w,y,a)$ is l.s.c.\ on $\X_W\times\A$ for each $y\in\X_Y.$

For an arbitrary fixed $y\in\X_Y$ the function
$(w,a)\mapsto \eta_{u}^\alpha(w,y,a)$ is \textit{$\K$-inf-compact on $\X_W\times\A$} as a sum of a $\K$-inf-compact function $(w,a)\mapsto c(w,y,a)$ and a nonnegative l.s.c.\ function $(w,a)\mapsto \alpha\int_\X  u(\tilde{x})q(d\tilde{x}|w,y,a)$ on $\X_W\times\A.$ Moreover, Berge's theorem for noncompact image sets \cite[Theorem~1.2]{FKN} implies that for each $(y,a)\in\X_Y\times\A$ the function $w\mapsto u^*(w,y):=\inf\limits_{a\in \A}\eta_u^\alpha(w,y;a)$ is l.s.c.\ on $\X_W.$
The Borel measurability of the function $u^*$ on $\X$ and the existence of a stationary policy $f:\X{\to} \A$ such that $u^*(x):=\eta_u^\alpha(x,f(x)),$ $x\in \mathbb{X},$ follow from \cite[Theorem~2.2 and Corollary~2.3(i)]{FK} because the function $(x,a)\mapsto \eta_{u}^\alpha(x,a)$ is Borel measurable on $\X\times\A$ and it is inf-compact in $a$ on $\A.$ Property~(a) for nonempty sets $\{A_*(x)\}_{x\in \X}$ follows from Borel measurability of $(x,a)\mapsto \eta_u^\alpha(x,a)$ on $\X\times\A$ and $x\mapsto u^*(x)$ on $\X.$ Property~(b) for  $\{A_*(x)\}_{x\in \X}$ follows from inf-compactness of $a\mapsto \eta_{u}^\alpha(x,a)$ on $\A$ for each $x\in\X.$
\end{proof}

\begin{proof}[Proof of Theorem~\ref{prop:dcoe}]
According to \cite[Proposition 8.2]{BS}, the functions $v_{t,\alpha}(x),$ $t=0,1,\ldots,$
recursively satisfy the optimality equations
with $v_{0,\alpha}(x)=0$ and
$v_{t+1,\alpha}(x)=\inf\limits_{a\in A(x)}\eta_{v_{t,\alpha}}^\alpha(x,a),$ for all $x\in
\X.$ So, Lemma~\ref{lm:NC1} sequentially applied to the functions $v_{0,\alpha}(x),$ $v_{1,\alpha}(x),\ldots,$ implies statement~(i) for them. According to \cite[Proposition~9.17]{BS}, $v_{t,\alpha}(x)\uparrow
v_{\alpha}(x)$ as $t \to +\infty$ for each $x\in\X.$ Therefore, $v_\alpha\in LW(\X).$ Thus, statement (i) is proved. In addition, \cite[Lemma 8.7]{BS} implies that a Markov policy defined at the first $T$ steps by the mappings $\phi_0^\alpha,...\phi_{T-1}^\alpha,$ that satisfy for all $t=1,\ldots,T$ the equations $v_{t,\alpha}(x)=\eta_{v_{t-1,\alpha}}^\alpha(x,\phi_{T-t}^\alpha(x)),$ for each $x\in \mathbb{X},$ is optimal for the horizon $T.$ According to \cite[Propositions~9.8 and 9.12]{BS}, $v_{\alpha}$ satisfies the discounted cost optimality equation $v_{\alpha}(x)=\inf\limits_{a\in A(x)}\eta_{v_{\alpha}}^\alpha(x,a)$ for each $x\in \X;$ and a stationary policy $\phi_\alpha$ is discount-optimal if and
only if $v_{\alpha}(x)=\eta_{v_{\alpha}}^\alpha(x,\phi_\alpha(x))$ for each $x\in \X.$  Statements  (ii-v) follow from these facts and Lemma~\ref{lm:NC1}.
\end{proof}

\begin{proof}[Proof of Theorem~\ref{th:mainMDPII}]
The equivalence of statements (a) and (b)  follows directly from Theorem~\ref{th:concept} applied to $\S_1:=\W,$ $\S_2:=\Y,$
$\S_3:=\W\times\A,$ and $\Psi:=P(\,\cdot\,|\,\cdot\,,y,\,\cdot\,).$
According to
\eqref{3.3}, Theorem~\ref{th:extra} applied to $\S_1:=\W,$ $\S_2:=\Y,$ $\S_3=\W,$ $\S_4:=\A,$ and $\Xi:=P(\,\cdot\,|\,\cdot\,,y,\,\cdot\,)$ implies that the stochastic kernel $P(\,\cdot\,|\,\cdot\,,y,\,\cdot\,)$ on $\W\times\Y$ given
$\W\times\A$ is {semi-uniform Feller} if and only if the stochastic kernel
$R(\,\cdot\,|\,\cdot\,,y,\,\cdot\,)$ on $\W\times\Y$ given $\P(\W)\times \A$ is {semi-uniform Feller}. Therefore,
statement~(b) holds if and only if the stochastic kernel $R(\,\cdot\,|\,\cdot\,,y,\,\cdot\,)$ on $\W\times\Y$ given
$\P(\W)\times\A$ is {semi-uniform Feller, that is, statement~(c) holds}. Thus, the equivalence of statements ({c})--({e}) 
follows directly from Theorem~\ref{th:CBRmain} applied to $\S_1:=\W,$ $\S_2:=\Y,$ $\S_3:=\P(\W)\times\A,$ $\Psi:=R(\,\cdot\,|\,\cdot\,,y,\,\cdot\,),$ $\Phi:=H(\,\cdot\,|\,\cdot\,,y,\,\cdot\,,\,\cdot\,),$ and $\phi:=q(\,\cdot\,|\,\cdot\,,y,\,\cdot\,).$

Moreover, let the nonnegative function $c$ be {$\M\K(\W\times\A,\Y)$}-inf-compact, and let for each $y\in\Y$ one of the equivalent conditions (a)--(d) hold. Then,  in view of \eqref{eq:c} and Theorem~\ref{th:PresC}, $\bar{c}$ is nonnegative and {$\M\K(\P(\W)\times\A,\Y)$}-inf-compact. Thus, the assumptions and  conclusions of Theorem~\ref{prop:dcoe} hold for the MDPCI $(\P(\W)\times\Y,\A,q,\c).$
\end{proof}

%

\begin{proof}[Proof of Theorem~\ref{th:mainMDPIINew}]
The equivalence of statements (a) and (b) follows directly from Theorem~\ref{th:concept} applied to $\S_1:=\W,$ $\S_2:=\Y,$ $\S_3:=\W\times\Y\times\A,$ and $\Psi:=P.$
According to \eqref{3.3}, Theorem~\ref{th:extra} applied to $\S_1:=\W,$ $\S_2:=\Y,$ $\S_3=\W,$ $\S_4:=\Y\times\A,$ and $\Xi:=P$ implies that the stochastic kernel $P$ on $\W\times\Y$ given
$\W\times\Y\times\A$ is {semi-uniform Feller} if and only if the stochastic kernel $R$ on $\W\times\Y$ given $\P(\W)\times\Y\times \A$ is {semi-uniform Feller}. Therefore,
statement~(b) holds if and only if the stochastic kernel $R$ on $\W\times\Y$ given $\P(\W)\times\Y\times \A$ is {semi-uniform Feller, that is, statement~(c) holds}. Thus, the equivalence of statements ({c})--({e}) 
{follows directly} from Theorem~\ref{th:CBRmain} applied to $\S_1:=\W,$ $\S_2:=\Y,$ $\S_3:=\P(\W)\times\Y\times\A,$ $\Psi:=R,$ $\Phi:=H,$ and $\phi:=q$.

Moreover, let the nonnegative function $c$ be {$\K$}-inf-compact, and let one of the equivalent conditions (a)--(d) hold. Then,  in view of \eqref{eq:c} and \cite[Theorem~3.3]{FKZ} on preservation of $\K$-inf-compactness, $\bar{c}$ is nonnegative and $\K$-inf-compact. Thus, according to Remark~\ref{rem:new}, the assumptions and  conclusions of Theorem~\ref{prop:dcoe} hold for the MDPCI $(\P(\W)\times\Y,\A,q,\c),$ and the functions $v_{t,\alpha},$ $t=0,1,\ldots,$ and $v_\alpha$ are lower semi-continuous.
\end{proof}

\begin{proof}[Proof of Corollary~\ref{cor:POMDPmain1}]
Let us prove that semi-uniform Feller continuity of the stochastic kernel $P$ on $\W\times\Y$ given $\W\times\A$ implies conditions (a) and (b). Indeed, Definition~\ref{defi:unifFP} implies weak continuity of the stochastic kernel $P_1$ on
$\W$ given $\W\times\A$ and continuity in the total variation of the stochastic kernel $Q_1$ on
$\Y$ given $\W\times\A$ because $P_1(\,\cdot\,|\,\cdot\,)=P(\,\cdot\,,\Y|\,\cdot\,)$ is weakly continuous and $Q_1(\,\cdot\,|\,\cdot\,)=P(\W,\,\cdot\,|\,\cdot\,)$ is continuous in total variation.  {Vice versa, let us prove} that conditions (a) and (b) imply semi-uniform Feller continuity of the stochastic kernel $P$ on $\W\times\Y$ given $\W\times\A.$
Indeed, $P$ on $\W\times\Y$ given
$\W\times\A$ is WTV-continuous since
\[
\begin{aligned}
&\ilim_{(w'a')\to(w,a)}\inf_{C\in\B(\Y)}(Q_1(C|w',a')P_1(\oo|w',a')-Q_1(C|w,a)P_1(\oo|w,a))\\
&\quad\ge \ilim_{(w'a')\to(w,a)}(P_1(\oo|w',a')-P_1(\oo|w,a))^-- \lim_{(w'a')\to(w,a)}\sup_{C\in\B(\Y)}|Q_1(C|w',a')-Q_1(C|w,a)|=0
\end{aligned}
\]
for each $\oo\in\tau(\W),$ where $a^-:=\min\{a,0\}$ for each $a\in\mathbb{R},$ the equality follows from weak continuity of $P_1$ on
$\W$ given $\W\times\A$ and  continuity in the total variation of $Q_1$ on
$\Y$ given $\A\times\W.$ Therefore, according to Theorem~\ref{th:equivWTV}(a,b), conditions (a) and (b) from Corollary~\ref{cor:POMDPmain1} taken together are equivalent to semi-uniform Feller continuity of the stochastic kernel $P$ on $\W\times\Y$ given $\W\times\A.$ Thus, Theorem~\ref{th:mainMDPII} {implies} all statements of Corollary~\ref{cor:POMDPmain1}.
\end{proof}

\begin{proof}[Proof of Corollary~\ref{cor:POMDPmain}]
For each $B\in \B(\W)$ consider the family of functions
\[
\mathcal{G}(B):=\Big\{(w,a)\mapsto \int_B Q_2(C|a,w')P_2(dw'|w,a)\,:\,C\in\B(\Y)\Big\}.
\]

Let condition~(a) hold. Fix an arbitrary open set $\oo\in\tau(\W).$ Feinberg et al. \cite[Theorem~1]{FKZJThProb}, applied to the lower semi-equicontinuous and uniformly bounded family of functions $\{(w',a)\mapsto\h\{w'\in\oo\}Q_2(C|a,w')\,:\,C\in\B(\Y)\}$ and weakly continuous  stochastic kernel $P_2(dw'|w,a)$ on $\W$ given $\W\times\A,$ implies that the family of functions $\mathcal{G}(\oo)$ is lower semi-equicontinuous at all the points $(w,a) \in \W \times \A,$ that is, the stochastic kernel $P$ on $\W\times\Y$ given $\W\times\Y\times\A$ defined in \eqref{eq:defi2} is WTV-continuous. Therefore, Theorem~\ref{th:equivWTV}(a,b) applied to the stochastic kernel $P$ on $\W\times\Y$ given $\W\times\Y\times\A$ implies that this kernel is semi-uniform Feller. Thus, assumption~(a) of Theorem~\ref{th:mainMDPII} holds, and this conclusion and Theorem~\ref{th:mainMDPII} imply all statements of Corollary~\ref{cor:POMDPmain} under condition~(a).


Now let condition~(b) hold. Let us prove that for each $B\in \B(\W)$ the family of functions
$\mathcal{G}(B)$ is equicontinuous at all $(w,a) \in \W \times \A,$ which implies condition~(a) of Theorem~\ref{th:mainMDPII}. Indeed,  for $n=1,2,\ldots,$
\begin{equation}\label{eq:Kara2a}
\begin{aligned}
\sup_{C\in\B(\Y)}\Big|\int_B Q_2(C|a^{\left(n\right)},w')P_2(dw'|w^{\left(n\right)},a^{\left(n\right)})-\int_B Q_2(C|a,w')P_2(dw'|w,a)\Big|\le I_1^{\left(n\right)}+I_2^{\left(n\right)},
\end{aligned}
\end{equation}
where $(w^{(n)},a^{(n)})\to (w,a)$ as $n\to\infty,$
\[
\begin{aligned}
&I_1^{\left(n\right)}:=\sup_{C\in\B(\Y)}\left|\int_B Q_2(C|a^{\left(n\right)},w')P_2(dw'|w^{\left(n\right)},a^{\left(n\right)})-\int_B Q_2(C|a^{\left(n\right)},w')P_2(dw'|w,a)\right|,\\
&I_2^{\left(n\right)}:=\sup_{C\in\B(\Y)}\int_B |Q_2(C|a^{\left(n\right)},w')-Q_2(C|a,w')|P_2(dw'|w,a).
\end{aligned}
\]
Let $C^{(n)}\in\B(\Y) $ be chosen to satisfy the inequality
\begin{equation}\label{eq:new1}
I_2^{\left(n\right)}\le
\int_B |Q_2(C^{\left(n\right)}|a^{\left(n\right)},w')-Q_2(C^{\left(n\right)}|a,w')|P_2(dw'|w,a)+\frac{1}{n},\qquad  n=1,2,\ldots.
\end{equation}
Note that $I_1^{\left(n\right)}\to 0$ as $n\to\infty$ because the {family} of measurable functions $\{w'\mapsto Q_2(C|\,a^{\left(n\right)},w')\,:\,n=1,2,\ldots\}$ is uniformly bounded by $1,$ and the stochastic kernel $P_2$ on $\W$ given $\W\times\A$ is continuous in total variation. Moreover, the convergence  $I_2^{\left(n\right)}\to 0$ as $n\to\infty$ follows from \eqref{eq:new1} and Lebesgue's dominated convergence theorem because the family of functions $\{w'\mapsto | Q_2(C^{\left(n\right)}|a^{\left(n\right)},w')-Q_2(C^{\left(n\right)}|a,w')|\,:\,n=1,2,\ldots\}$ is uniformly bounded by $1$  and pointwise convergent to $0,$ according to \eqref{eq:Kara1}. Therefore, the family of functions $\mathcal{G}(B)$ is equicontinuous on $\W \times \A.$ Thus, assumption~(a) of Theorem~\ref{th:mainMDPII} holds, and this conclusion and Theorem~\ref{th:mainMDPII} imply all statements of Corollary~\ref{cor:POMDPmain} under condition~(b).
\end{proof}

\section*{Acknowledgements} We thank Janey (Huizhen) Yu for valuable remarks. Research of the second and the third authors was partially supported by the National Research Foundation of Ukraine, Grant No.~2020.01/0283.  We thank the referees for insightful remarks. In particular, one of the referees suggested a short proof of weakly continuity of semi-uniform Feller kernels, observed the equivalence of WTV-continuity and semi-uniform Feller continuity,  suggested to strengthen Theorems~\ref{prop:dcoe} and \ref{th:mainMDPII} to their current formulations, proposed the provided proof of Theorem~\ref{th:CBRmain}, and made other valuable comments.

\end{document}